\documentclass[reqno,12pt]{amsart}
\usepackage{graphicx}
\usepackage{epsfig}
\usepackage{url}
\usepackage{amssymb,amsthm,amsmath,amsfonts,amstext,mathrsfs}
\usepackage{latexsym}
\usepackage[all]{xy}
\newtheorem{thm}{Theorem}[section]
\newtheorem{cor}[thm]{Corollary}
\newtheorem{lem}[thm]{Lemma}
\newtheorem{prop}[thm]{Proposition}
\newtheorem{exper}[thm]{Experimental observation}
\newtheorem{defin}[thm]{Definition}

\input xy
\xyoption{all}

\def\d{\,{\rm{d}}}
\def\p{\,{\partial}}
\def\pp{\,{{\sf p}}}

\def\fp{\,{F_{{\sf p}}}}
\def\I{\,{\mathscr{I}}}

\title[Explicit series for the dyadic period function]
{The Minkowski question mark function: explicit series for the dyadic period function and moments}
\author[Giedrius Alkauskas]{Giedrius Alkauskas}
\begin{document}

\begin{abstract} Previously, several natural integral transforms of
the Minkowski question mark function $F(x)$ were introduced by the author. Each of
them is uniquely characterized by certain regularity conditions and the
functional equation, thus encoding intrinsic information about $F(x)$. One of
them - the dyadic period function $G(z)$ - was defined as a Stieltjes transform. In this paper we introduce a family of
``distributions" $\fp(x)$ for $\Re\pp\geq 1$, such that $F_{1}(x)$ is the
question mark function and $F_{2}(x)$ is a discrete distribution with support
on $x=1$. We prove that the generating function of moments of $F_{\pp}(x)$ satisfies the three term functional equation.
This has an independent interest, though our main concern is the information it provides about $F(x)$. This approach yields
the following main result: we prove that the dyadic period function is a sum of infinite series of rational
functions with rational coefficients.
\end{abstract}
\maketitle
\begin{center}
\end{center}
\begin{center}
\end{center}
\setcounter{tocdepth}{2} \tableofcontents
\textbf{Acknowledgements.} The author sincerely thanks J\"{o}rn Steuding, whose
seemingly elementary problem, proposed at the problem session in Palanga
conference in 2006, turned to be a very rich and generous one. Also, the author thanks several colleagues for showing interest in these results, after
preprints became available in {\tt arXiv}, especially to Steven Finch, Jeffrey Lagarias and Stefano Isola.
\section{Introduction and main result}
The aim of this paper to continue investigations on the moments of Minkowski
$?(x)$ function, begun in \cite{ga1}, \cite{ga2} and \cite{ga3}. The function
$?(x)$ (``the question mark function") was introduced by Minkowski as an
example of a continuous function $F:[0,\infty)\rightarrow[0,1)$, which maps
rationals to dyadic rationals, and quadratic irrationals to non-dyadic
rationals. For non-negative real $x$ it is defined by the expression
\begin{eqnarray}
F([a_{0},a_{1},a_{2},a_{3},...])=1-2^{-a_{0}}+2^{-(a_{0}+a_{1})}-2^{-(a_{0}+a_{1}+a_{2})}+...,
\label{min}
\end{eqnarray}
where $x=[a_{0},a_{1},a_{2},a_{3},...]$ stands for the representation of $x$ by
a (regular) continued fraction \cite{khin}. By tradition, this function is more
often investigated in the interval $[0,1]$, and in this case it is normalized
in order $F(1)=1$, whereas in our case $F(1)=\frac{1}{2}$. Accordingly, we make
a convention that $?(x)=2F(x)$ for $x\in[0,1]$. For rational $x$, the series
terminates at the last nonzero partial quotient $a_{n}$ of the continued fraction. This
function is continuous, monotone and singular \cite{denjoy}. By far not
complete overview of the papers written about the Minkowski question mark
function or closely related topics (Farey tree, enumeration of rationals,
Stern's diatomic sequence, various 1-dimensional generalizations and
generalizations to higher dimensions, statistics of denominators and Farey
intervals, Hausdorff dimension and analytic properties) can be found in
\cite{ga1}. These works include \cite{bonano_graffi_isola}, \cite{bonnano}, \cite{cw}, \cite{denjoy},
\cite{dm}, \cite{girgensohn}, \cite{grabner_tichy} (this is the only paper where the moments of a certain singular distribution -
a close relative of $F(x)$ - were considered), \cite{isola}, \cite{kess},
\cite{kinney}, \cite{lag}, \cite{lagtresser}, \cite{okamoto_wunsch}, \cite{pan}
\cite{paradis1}, \cite{paradis2}, \cite{ramharter}, \cite{reese}, \cite{ryde2}, \cite{salem}, \cite{tichy_uitz}.
The internet page \cite{mink} contains up-to-date and exhaustive bibliography list of
papers related to Minkowski question mark function.\\
\indent Recently, in Calkin and Wilf \cite{cw} defined a binary tree which is
generated by the iteration
$$
{a\over b}\quad\mapsto\quad {a\over a+b}\ ,\quad {a+b\over b},
$$
starting from the root ${1\over 1}$. The last two authors have greatly publicized this tree, but it was known long ago to physicists and mathematicians (alias, Stern-Brocot or Farey tree). Elementary considerations show that this
tree contains every positive rational number once and only once, each being
represented in lowest terms. The first four iterations lead to
\begin{eqnarray}
\xymatrix @R=.5pc @C=.5pc { & & & & & & & {1\over 1} & & & & & & & \\
& & & {1\over 2} \ar@{-}[urrrr] & & & & & & & & {2\over 1} \ar@{-}[ullll] & & & \\
& {1\over 3} \ar@{-}[urr] & & & & {3\over 2}\ar@{-}[ull] & & & & {2\over 3}\ar@{-}[urr] & & & & {3\over 1} \ar@{-}[ull] & \\
{1\over 4} \ar@{-}[ur] & & {4\over 3} \ar@{-}[ul] & & {3\over 5}
\ar@{-}[ur] & & {5\over 2} \ar@{-}[ul] & & {2\over 5} \ar@{-}[ur]
& & {5\over 3} \ar@{-}[ul] & & {3\over 4} \ar@{-}[ur] & & {4\over
1} \ar@{-}[ul] }
\label{cw}
\end{eqnarray}
\indent It is of utmost importance to note that the $n$th generation consists
of exactly those $2^{n-1}$ positive rational numbers, whose elements of the
continued fraction sum up to $n$. This fact can be easily inherited directly
from the definition. First, if rational number $\frac{a}{b}$ is represented as
a continued fraction $[a_{0},a_{1},...,a_{r}]$, then the map
$\frac{a}{b}\rightarrow\frac{a+b}{b}$ maps $\frac{a}{b}$ to
$[a_{0}+1,a_{1}...,a_{r}]$. Second, the map
$\frac{a}{b}\rightarrow\frac{a}{a+b}$ maps $\frac{a}{b}$ to
$[0,a_{1}+1,...,a_{r}]$ in case $\frac{a}{b}<1$, and to
$[1,a_{0},a_{1},...,a_{r}]$ in case $\frac{a}{b}>1$. This is an important fact
which makes the investigations of rational numbers according to their position
in the Calkin-Wilf tree highly motivated from the perspective of metric
number theory and dynamics of continued fractions. \\
\indent It is well known that each generation of (\ref{cw}) possesses a distribution
function $F_{n}(x)$, and $F_{n}(x)$ converges uniformly to $F(x)$. The function
$F(x)$ as a distribution function (in the sense of probability theory, which
imposes the condition of monotonicity) is uniquely determined by the functional
equation \cite{ga1}
\begin{eqnarray}
2F(x)=\left\{\begin{array}{c@{\qquad}l} F(x-1)+1 & \mbox{if}\quad
x\geq 1,
\\ F({x\over 1-x}) & \mbox{if}\quad 0\leq x<1. \end{array}\right.
\label{distr}
\end{eqnarray}
This implies $F(x)+F(1/x)=1$. The mean value of
$F(x)$ has been investigated by several authors, and was proved to be $3/2$.\\
\indent Lastly, and most importantly, let us point out that, surprisingly,
there are striking similarities and parallels between the results proved in
\cite{ga1} and \cite{ga2} with Lewis'-Zagier's (\cite{lewis}, \cite{zagl})
results on period functions for Maass wave forms. (see \cite{ga2} for the
explanation of this
phenomena).\\

Just before formulating the main Theorem of this paper, we provide a short
summary of previous results proved by the author about certain natural integral
transforms of $F(x)$. Let
\begin{eqnarray*}
M_{L}=\int\limits_{0}^{\infty}x^{L}\d F(x),\quad
m_{L}=\int\limits_{0}^{\infty}\Big{(}\frac{x}{x+1}\Big{)}^{L}\d
F(x)=2\int\limits_{0}^{1}x^{L}\d F(x).
\end{eqnarray*}
Both sequences are of definite number-theoretical significance because
\begin{eqnarray*}
M_{L}=\lim_{n\rightarrow\infty}2^{1-n}\sum\limits_{a_{0}+a_{1}+...+a_{s}=n}[a_{0},a_{1},..,a_{s}]^{L},\quad
m_{L}=\lim_{n\rightarrow\infty}2^{2-n}\sum\limits_{a_{1}+...+a_{s}=n}[0,a_{1},..,a_{s}]^{L},
\end{eqnarray*}
(the summation takes place over rational numbers represented as continued
fractions; thus, $a_{i}\geq 1$ and $a_{s}\geq 2$). We define the exponential
generating functions
\begin{eqnarray*}
M(t)&=&\sum_{L=0}^{\infty}\frac{M_{L}}{L!}t^{L}=\int\limits_{0}^{\infty}e^{xt}\d
F(x),\\
\mathfrak{m}(t)&=&\sum_{L=0}^{\infty}\frac{m_{L}}{L!}t^{L}=\int\limits_{0}^{\infty}\exp\Big{(}\frac{xt}{x+1}\Big{)}\d
F(x)=2\int\limits_{0}^{1}e^{xt}\d F(x).
\end{eqnarray*}
One directly verifies that $\mathfrak{m}(t)$ is an entire function, and that
$M(t)$ is meromorphic function with simple poles at $z=\log 2+2\pi i n$,
$n\in\mathbb{Z}$. Further, we have
\begin{eqnarray*}
M(t)=\frac{\mathfrak{m}(t)}{2-e^{t}},\quad
\mathfrak{m}(t)=e^{t}\mathfrak{m}(-t).
\end{eqnarray*}
 The second identity represents only the symmetry property, given by
$F(x)+F(1/x)=1$. The main result about $\mathfrak{m}(t)$ is that it is uniquely
determined by the regularity condition $\mathfrak{m}(-t)\ll e^{-\sqrt{\log
2}\sqrt{t}}$, as $t\rightarrow\infty$, the boundary condition $\mathfrak{m}(0)=1$, and the
integral equation
\begin{eqnarray}
\mathfrak{m}(-s)=(2e^{s}-1)\int\limits_{0}^{\infty}\mathfrak{m}'(-t)J_{0}(2\sqrt{st})\d
t,\quad s\in\mathbb{R}_{+}.
\label{char}
\end{eqnarray}
(Here $J_{0}(\star)$ stands for the Bessel function
$J_{0}(z)=\frac{1}{\pi}\int_{0}^{\pi}\cos(z\sin x)\d x$).\\

Our primary object of investigations is the generating function of moments. Let
$G(z)=\sum\limits_{L=1}^{\infty}m_{L}z^{L-1}$. This series converges for
$|z|\leq 1$, and the functional equation for $G(z)$ (see below) implies that
there exist all derivatives of $G(z)$ at $z=1$, if we approach this point while
remaining in the domain $\Re z\leq 1$. Then the integral
\begin{eqnarray}
G(z)=\int\limits_{0}^{\infty}\frac{1}{x+1-z}\d
F(x)=2\int\limits_{0}^{1}\frac{x}{1-xz}\d F(x). \label{gen}
\end{eqnarray}
(which is Stieltjes transform  of $F(x)$) extends $G(z)$ to the cut plane $\mathbb{C}\setminus(1,\infty)$. The generating
function of moments $M_{L}$ does not exist due to the factorial growth of
$M_{L}$, but this generating function can still be defined in the cut plane
$\mathbb{C}'=\mathbb{C}\setminus(0,\infty)$ by
$\int_{0}^{\infty}\frac{x}{1-xz}\d F(x)$. In fact, this integral just equals to
$G(z+1)$. Thus, there exist all higher derivatives of $G(z)$ at $z=1$, and
$\frac{1}{(L-1)!}\frac{\d^{L-1}}{\d z^{L-1}}G(z)\big{|}_{z=1}=M_{L}$, $L\geq 1$.
The following result was proved in \cite{ga1}.
\begin{thm} The function $G(z)$, defined initially as a power series, has an analytic
continuation to the cut plane $\mathbb{C}\setminus(1,\infty)$ via (\ref{gen}).
It satisfies the functional equation
\begin{eqnarray}
\frac{1}{z}+\frac{1}{z^{2}}G\Big{(}\frac{1}{z}\Big{)}+2G(z+1)=G(z),\label{funk}
\end{eqnarray}
and also the symmetry property
\begin{eqnarray*}
G(z+1)=-\frac{1}{z^{2}}G\Big{(}\frac{1}{z}+1\Big{)}-\frac{1}{z},
\end{eqnarray*}
(which is a consequence of the main functional equation). Moreover,
$G(z)\rightarrow 0$, if $z\rightarrow\infty$ and the distance from $z$ to a
half line $[0,\infty)$ tends to infinity.
Conversely, the function having these properties is unique.
\label{thm1.0}
\end{thm}
\indent Accordingly, this result and the specific appearance of the three term
functional equation justifies the name for $G(z)$ as {\it the dyadic period
function}.\\

We wish to emphasize that the main motivation for previous research was
clarification of the nature and structure of the moments $m_{L}$. It was
greatly desirable to give these constants (emerging as if from geometric chaos)
some other expression than the one obtained directly from the Farey (or
Calkin-Wilf) tree, which could reveal their structure to greater extent. This
is accomplished in the current work. Thus, the main result can be formulated as
follows.
\begin{thm}
There exist canonical and explicit sequence of rational functions
$\mathbf{H}_{n}(z)$, such that for $\{|z|\leq \frac{3}{4}\}\cup\{|z+\frac{9}{7}|\leq \frac{12}{7}\}$ one has an absolutely convergent series
\begin{eqnarray*}
G(z)=\int\limits_{0}^{\infty}\frac{1}{x+1-z}\d
F(x)=\sum\limits_{n=0}^{\infty}(-1)^{n}\mathbf{H}_{n}(z),\quad
\mathbf{H}_{n}(z)=\frac{\mathscr{B}_{n}(z)}{(z-2)^{n+1}},
\end{eqnarray*}
where $\mathscr{B}_{n}(z)$ is polynomial with rational coefficients of degree
$n-1$. For $n\geq 1$ it has the following reciprocity property:
\begin{eqnarray*}
\mathscr{B}_{n}(z+1)=(-1)^{n}z^{n-1}\mathscr{B}_{n}\Big{(}\frac{1}{z}+1\Big{)},\quad\mathscr{B}_{n}(0)=0.
\end{eqnarray*}
\label{thm1.1}
\end{thm}
The rational function $\mathbf{H}_{n}(z)$ are defined via implicit and rather complicated recurrence (\ref{lygtis}) (see Section 6).
The following table gives initial polynomials $\mathscr{B}_{n}(z)$. \\

\noindent\begin{tabular}{|r | r || r| r|} \hline
$n$   & $\mathscr{B}_{n}(z)$                    & $n$ & $\mathscr{B}_{n}(z)$\\
\hline
& & &\\
$0$&$-1$            &$4$ & $\displaystyle-\frac{2}{27}z^{3}+\frac{53}{270}z^{2}-\frac{53}{270}z$\\
& & &\\
$1$&$0$            &$5$&$\displaystyle\frac{4}{81}z^{4}-\frac{104}{675}z^{3}+\frac{112}{675}z^{2}-\frac{224}{2025}z$\\
&  & &\\
$2$&$\displaystyle-\frac{1}{6}z$&$6$&$\displaystyle-\frac{8}{243}z^{5}+\frac{47029}{425250}z^{4}-\frac{1384}{14175}z^{3}-\frac{787}{30375}z^{2}+\frac{787}{60750}z$\\
& & & \\
$3$&$\displaystyle\frac{1}{9}z^{2}-\frac{2}{9}z$&$7$&
$\displaystyle\frac{16}{729}z^6-\frac{1628392}{22325625}z^5+\frac{272869}{22325625}z^4+
\frac{5392444}{22325625}z^3-\frac{238901}{637875}z^2+\frac{477802}{3189375}z$\\
& & &\\
\hline
\end{tabular}\\

\noindent {\it Remark. }The constant $\frac{3}{4}$ can be replaced by any constant less than $1.29^{-1}$ (the latter comes exactly from Lemma \ref{lema.3}).
Unfortunately, our method does not allow to prove an absolute convergence in the disk $|z|\leq 1$.
In fact, apparently the true region of convergence of the series
in question is the half plane $\Re z\leq 1$. Take, for example,
$z_{0}=\frac{2}{3}+4i$. Then by (\ref{funk}) and symmetry property one has
\begin{eqnarray*}
G(z_{0})=\frac{1}{2}G(z_{0}-1)-\frac{1}{2(z_{0}-1)^{2}}G\Big{(}\frac{1}{z_{0}-1}\Big{)}-\frac{1}{2(z_{0}-1)}=\\
-\frac{1}{2(z_{0}-2)^2}G\Big{(}\frac{z_{0}-1}{z_{0}-2}\Big{)}-
\frac{1}{2(z_{0}-1)^{2}}G\Big{(}\frac{1}{z_{0}-1}\Big{)}-\frac{1}{2(z_{0}-2)}-\frac{1}{2(z_{0}-1)}.
\end{eqnarray*}
Both arguments under $G$ on the right belong to the unit circle, and thus we
can use Taylor series for $G(z)$. Using numerical values of $m_{L}$, obtained
via the method described in Appendix A.2., we obtain:
$G(z_{0})=0.078083_{+}+0.205424_{+}i$, with all digits exact. On the other
hand, the series in Theorem \ref{thm1.1} for $n=60$ gives
\begin{eqnarray*}
\sum\limits_{n=0}^{60}(-1)^{n}\mathbf{H}_{n}(z_{0})=0.078090_{+}+0.205427_{+}i.
\end{eqnarray*}
Finally, based on the last integral in (\ref{gen}), we can calculate $G(z)$ as
a Stieltjes integral. If we divide the unit interval into $N=3560$ equal
subintervals, and use Riemann-Stieltjes sum, we get an approximate value
$G(z_{0})\approx 0.078082_{+}+0.205424_{+}i$. All evaluations match very well.
\begin{exper}
We conjecture that the series in Theorem \ref{thm1.1} converges absolutely for $\Re z\leq 1$.
\label{exper1.2}
\end{exper}
 With a slight abuse of notation, we will henceforth write
$f^{(L-1)}(z_{0})$ instead of $\frac{\d^{L-1}}{\d
z^{L-1}}f(z)\big{|}_{z=z_{0}}$.
\begin{cor}
The moments $m_{L}$ can be expressed by the convergent
series of rational numbers:
\begin{eqnarray*}
m_{L}=\lim_{n\rightarrow\infty}2^{2-n}\sum\limits_{a_{1}+a_{2}+...+a_{s}=n}[0,a_{1},a_{2},...,a_{s}]^{L}=
\frac{1}{(L-1)!}\sum\limits_{n=0}^{\infty}(-1)^{n}\mathbf{H}^{(L-1)}_{n}(0).
\end{eqnarray*}
The speed of convergence is given by the following estimate:
$\Big{|}\mathbf{H}^{(L-1)}_{n}(0)\Big{|}\ll\frac{1}{n^{M}}$, for every
$M\in\mathbb{N}$. The implied constant depends only on $L$ and $M$.
\label{cor1.3}
\end{cor}
Thus, $m_{2}=\sum_{n=0}^{\infty}(-1)^{n}\mathbf{H}'_{n}(0)=0.2909264764_{+}$.
Regarding the speed, numerical calculations show that in fact the convergence
is geometric. Theorem \ref{thm1.1} in case $z=1$ gives
\begin{eqnarray*}
M_{1}=G(1)=1+0+\sum\limits_{n=0}^{\infty}\frac{1}{6}\Big{(}\frac{2}{3}\Big{)}^{n}=\frac{3}{2},
\end{eqnarray*}
which we already know (see Corollary \ref{cor4.5}; the above is a Taylor series for $M_{1}(\pp)$ in powers of $\pp-2$,
specialized at $\pp_{0}=1$).
Geometric convergence would be the consequence of the
fact that analytic functions $m_{L}(\pp)$ extend beyond $\pp=1$ (see below).
This is supported by the phenomena represented as Experimental observation \ref{exper1.4}. Meanwhile, we are
able to prove only the given rate. If we were allowed to use the point $z=1$, Theorem \ref{thm1.1} would
give a convergent series for the
moments $M_{L}$ as well. This is exactly the same as the series in the
Corollary \ref{cor1.3}, only one needs to use a point $z=1$ instead of $z=0$.
\begin{exper}
For $L\geq 1$, the series
\begin{eqnarray*}
M_{L}(\pp)=\frac{1}{(L-1)!}\cdot\sum\limits_{n=0}^{\infty}(\pp-2)^{n}\mathbf{H}^{(L-1)}_{n}(1),\quad M_{L}(1)=M_{L},
\end{eqnarray*}
has exactly $2-\frac{1}{\sqrt[L]{2}}$ as a radius of convergence.
\label{exper1.4}
\end{exper}
To this account, Proposition \ref{prop4.3} endorse this phenomena, which is
highly supported by numerical calculations, and which does hold for $L=1$.\\

The following two tables give starting values for the sequence $\mathbf{H}'_{n}(0)$.\\

\noindent\begin{center}
\begin{tabular}{|r | c||r|c||r| c|}
\hline
$n$ & $\mathbf{H}'_{n}(0)$& $n$&$\mathbf{H}'_{n}(0)$& $n$&$\mathbf{H}'_{n}(0)$\\
\hline
& & & & &\\
$0$ & $\displaystyle\frac{1}{4}$    &$5$&$\displaystyle-\frac{7}{2\cdot 3^4\cdot 5^2}$                  &$10$&$\displaystyle-\frac{8026531718888633}{2^{12}\cdot 3^{9}\cdot 5^{7}\cdot 7^{4}\cdot 11\cdot  17^{2}}$\\
& & & & &\\
$1$ & $\displaystyle0$              &$6$&$\displaystyle-\frac{787}{2^8\cdot 3^5\cdot 5^3}$              &$11$&$\displaystyle\frac{797209536976557079423}{2^{11}\cdot 3^{10}\cdot 5^{8}\cdot 7^{5}\cdot 11^{2}\cdot 17^{3}\cdot 31}$\\
& & & & &\\
$2$ & $\displaystyle\frac{1}{48}$   &$7$&$\displaystyle\frac{238901}{2^7\cdot 3^6\cdot 5^4\cdot 7}$          &$12$&$\displaystyle\frac{4198988799919158293319845971}{2^{14}\cdot 3^{11}\cdot 5^{9}\cdot 7^{6}\cdot 11^{3}\cdot 13\cdot  17^{4}\cdot 31^{2}}$\\
& & & & &\\
$3$ & $\displaystyle-\frac{1}{72}$  &$8$&$\displaystyle-\frac{181993843}{2^{10}\cdot 3^7\cdot 5^5\cdot 7^2}$   &$13$&$\displaystyle-\frac{12702956822417247965298252330349561}{2^{10}\cdot 3^{12}\cdot 5^{10}\cdot 7^{7}\cdot 11^{4}\cdot 13^{2}\cdot 17^{5}\cdot 31^{3}}$\\
& & & & &\\
$4$&$\displaystyle\frac{53}{8640}$  &$9$&$\displaystyle\frac{12965510861}{2^6\cdot 3^8\cdot 5^6\cdot 7^3\cdot 17}$&$14$&$\displaystyle\frac{7226191636013675292833514548603516395499899}{2^{16}\cdot 3^{13}\cdot 5^{11}\cdot 7^{8}\cdot 11^{5}\cdot 13^{3}\cdot 17^{6}\cdot 31^{4}}$\\
& & & & &\\
\hline
\end{tabular}
\end{center}
\bigskip

\noindent\begin{center}
\begin{tabular}{|r | c|}
\hline
$n$ & $\mathbf{H}'_{n}(0)$\\
\hline
&\\
$15$ &
$-\displaystyle\frac{129337183009042141853748450730581369733226857443915617}{2^{15}\cdot
3^{14}\cdot 5^{12}\cdot 7^{9}\cdot 11^{6}\cdot 13^{4}\cdot 17^{7}\cdot
31^{5}\cdot 43\cdot 127}$\\
&\\
$16$ &
$\displaystyle\frac{31258186275777197041073243752715109842753785598306812028984213251}{2^{18}\cdot
3^{15}\cdot 5^{13}\cdot 7^{10}\cdot 11^{7}\cdot 13^{5}\cdot 17^{8}\cdot
31^{6}\cdot 43^{2}\cdot 127^{2}}$\\
&\\
$17$ &
$-\displaystyle\frac{3282520501229639755997762022707321704397776888948469860959830459774414444483}{2^{12}\cdot
3^{16}\cdot 5^{14}\cdot 7^{11}\cdot 11^{8}\cdot 13^{6}\cdot 17^{9}\cdot
31^{7}\cdot 43^{3}\cdot 127^{3}\cdot 257}$\\
&\\
\hline
\end{tabular}
\end{center}
\bigskip
\indent The float values of the last three rational numbers are $-0.000025804822076$, $0.000018040274062$ and $-0.000010917558446$ respectively. The alternating sum of the elements in the table is
$\sum_{n=0}^{N}(-1)^{n}\mathbf{H}'_{n}(0)=0.2909255862_{+}$ (where $N=17$),
whereas $N=40$ gives $0.2909264880_{+}$, and $N=50$ gives
$0.2909264784_{+}$. Note that the manifestation of Fermat and Mersenne primes in the denominators at an early stage is not accidental,
minding the exact value of the determinant in Lemma \ref{lem6.1}, Chapter 6 (see below).
Moreover, the prime powers of every odd prime, which divides the denominator, increase each time by $1$ while
passing from $\mathbf{H}'_{n}(0)$ to $\mathbf{H}'_{n+1}(0)$. The pattern for the powers of $2$ is more complicated. More thorough research
of the linear map in Lemma \ref{lem6.1} can thus clarify prime decomposition of denominators; numerators remains much more complicated.
\\

As will be apparent later, the result in Theorem \ref{thm1.1} is derived from the
knowledge of $\pp-$derivatives of $G(\pp,z)$ at $\pp=2$ (see below). On the
other hand, since there are two points ($\pp=2$ and $\pp=0$) such that all
higher $\pp-$derivatives of $G(\pp,z)$ are rational functions in $z$, it is not
completely surprising that the approach through $\pp=0$ also gives convergent
series for the moments, though in this case we are forced to use Borel
summation. At this point, the author does not have a strict mathematical proof
of this result (since the function $G(\pp,z)$ is meanwhile defined only for
$\Re \pp\geq 1$), though numerical calculations provide overwhelming evidence
for its validity.
\begin{exper}Define the rational functions (with rational coefficients) $\mathbf{Q}_{n}(z)$, $n\geq 0$, by
\begin{eqnarray*}
\mathbf{Q}_{0}(z)=-\frac{1}{2z},\text{ and recurrently by }
\mathbf{Q}_{n}(z)=\frac{1}{2}\sum\limits_{j=0}^{n-1}\frac{1}{j!}\cdot\frac{\p^{j}}{\p
z^{j}}\mathbf{Q}_{n-j-1}(-1)\cdot\Big{(}z^{j}-\frac{1}{z^{j+2}}\Big{)}.
\end{eqnarray*}
Then
\begin{eqnarray}
m_{L}=\lim_{n\rightarrow\infty}2^{2-n}\sum\limits_{a_{1}+a_{2}+...+a_{s}=n}[0,a_{1},a_{2},...,a_{s}]^{L}=
\frac{1}{(L-1)!}
\sum\limits_{r=0}^{\infty}\Big{(}\sum\limits_{n=0}^{\infty}\frac{\mathbf{Q}^{(L-1)}_{n}(-1)}{n!}\cdot\int\limits_{r}^{r+1}t^{n}e^{-t}\d
t\Big{)}.\label{exr}
\end{eqnarray}
Moreover,
\begin{eqnarray*}
\mathbf{Q}_{n}(z)=\frac{(z+1)(z-1)\mathscr{D}_{n}(z)}{z^{n+1}},\quad n\geq 1,
\end{eqnarray*}
where $\mathscr{D}_{n}(z)$ are polynomials with rational coefficients
($\mathbb{Q}_{p}$ integers for $p\neq 2$) of degree $2n-2$ with the reciprocity
property
\begin{eqnarray*}
\mathscr{D}_{n}(z)=z^{2n-2}\mathscr{D}_{n}\Big{(}\frac{1}{z}\Big{)}.
\end{eqnarray*}
\label{exper1.5}
\end{exper}
Note the order of summation in the series for $m_{L}$, since the reason for
introducing exponential function is because we use Borel summation. For
example,
\begin{eqnarray*}
``1-2+4-8+16-32+..."\mathop{=}^{{\tt Borel}}\sum\limits_{r=0}^{\infty}\Big{(}\sum\limits_{n=0}^{\infty
}\frac{(-2)^{n}}{n!}\cdot\int\limits_{r}^{r+1}t^{n}e^{-t}\d
t\Big{)}=\frac{1}{3}.
\end{eqnarray*}
The following table gives initial results.\\

\begin{center}
\begin{tabular}{|r | l || r| l|}
\hline
$n$   & $\mathscr{D}_{n}(z)$                    & $n$ & $\mathscr{D}_{n}(z)$\\
\hline
$1$      & $\frac{1}{4}$                          & $4$ & $\frac{1}{8}(2z^{6}-3z^{5}+6z^{4}-3z^{3}+6z^{2}-3z+2)$\\
$2$      & $\frac{1}{4}(z^{2}+1)$               & $5$ & $\frac{1}{4}(z^8-2z^{7}+4z^{6}-7z^{5}+4z^{4}-7z^{3}+4z^{2}-2z+1)$\\
$3$    & $\frac{1}{4}(z^{4}-z^{3}+z^{2}-z+1)$               & $6$  & $\frac{1}{8}(2z^{10}-5z^{9}+12z^{8}-20z^{7}+37z^{6}-$\\
   &                                           &     & $-20z^{5}+37z^{4}-20z^{3}+12z^{2}-5z+2)$\\
\hline
\end{tabular}\\
\end{center}
\bigskip
 The next table gives
$\mathbf{Q}'_{n}(-1)=2(-1)^{n}\mathscr{D}_{n}(-1)$ explicitly: these constants appear in
the series defining the first non-trivial moment $m_{2}$. Also, since these
numbers are $p-$adic integers for $p\neq 2$, there is a hope for the successful
implementation of the idea from the last chapter in \cite{ga2}; that
is, possibly one can define moments $m_{L}$ as $p-$adic rationals as well.\\

\begin{center}
\begin{tabular}{|r r||r r||r r||r c|}
\hline
$n$ & $\mathbf{Q}'_{n}(-1)$& $n$&$\mathbf{Q}'_{n}(-1)$& $n$&$\mathbf{Q}'_{n}(-1)$& $n$ & $\mathbf{Q}'_{n}(-1)$\\
\hline
$0$ & $\frac{1}{2}$       &$8$ & $\frac{1417}{4}$        &$16$& $\frac{206836175}{64}$     &$24$&  $\frac{1685121707817}{32}$            \\
$1$ & $-\frac{1}{2}$      &$9$ & $-\frac{8431}{8}$       &$17$&$-\frac{339942899}{32}$     &$25$&  $-\frac{92779913448103}{512}$          \\
$2$ & $1$                 &$10$& $\frac{50899}{16}$      &$18$&$\frac{1125752909}{32}$     &$26$&  $\frac{80142274019997}{128}$          \\
$3$ & $-\frac{5}{2}$      &$11$& $-9751$                 &$19$&$-\frac{15014220659}{128}$  &$27$&  $-\frac{1111839248032133}{512}$          \\
$4$ & $\frac{25}{4}$      &$12$& $30365$                 &$20$&$\frac{25188552721}{64}$    &$28$&  $\frac{7740056893342455}{1024}$          \\
$5$ & $-16$               &$13$& $-\frac{3069719}{32}$   &$21$&$-\frac{170016460947}{128}$ &$29$&  $-\frac{13515970598654393}{512}$          \\
$6$ & $43$                &$14$& $\frac{1227099}{4}$     &$22$&$\frac{1153784184807}{256}$ &$30$&  $\frac{47354245650630005}{512}$          \\
$7$ & $-\frac{971}{8}$    &$15$& $-\frac{31719165}{32}$  &$23$&$-\frac{983668214037}{64}$  &$31$&  $-\frac{665632101181145115}{2048}$       \\
\hline
\end{tabular}\\
\end{center}
\bigskip
The final table in this section lists float values of the constants
\begin{eqnarray*}
\vartheta_{r}=\sum\limits_{n=0}^{\infty}\frac{\mathbf{Q}'_{n}(-1)}{n!}
\cdot\int\limits_{r}^{r+1}t^{n}e^{-t}\d t,\quad r\in\mathbb{N}_{0},\quad
\sum\limits_{r=0}^{\infty}\vartheta_{r}=m_{2},
\end{eqnarray*}
appearing in Borel summation.
\begin{center}
\begin{tabular}{|r | r || r| r|}
\hline
$r$   & $\vartheta_{r}$                    & $r$ & $\vartheta_{r}$\\
\hline
$0$ & $0.2327797875$  &  $6$ & $0.0004701146$\\
$1$ & $0.0471561089$  &  $7$ & $0.0004980015$\\
$2$ & $0.0085133626$  &  $8$ & $0.0004005270$\\
$3$ & $0.0005892453$  &  $9$ & $0.0002722002$\\
$4$ & $-0.0001872357$ & $10$ & $0.0001607897$\\
$5$ & $0.0002058729$  & $11$ & $0.0000812407$\\
\hline
\end{tabular}
\end{center}
\bigskip
Thus, $\sum_{r=0}^{11}\vartheta_{r}=0.2909400155_{+}=m_{2}+0.000013539_{+}$.\\

This paper is organized as follows. In Section 2, for each $\pp$, $1\leq
\pp<\infty$, we introduce a generalization of the Farey (Calkin-Wilf) tree, denoted by
$\mathcal{Q}_{\pp}$. This leads to the notion of $\pp-$continued fractions and
$\pp-$Minkowski question mark functions $\fp(x)$. Though $\pp-$continued
fractions are of independent interest (one could define a transfer operator, to
prove an analogue of Gauss-Kuzmin-L\'{e}vy theorem, various metric results and
introduce structural constants), we confine to the facts which are necessary
for our purposes and leave the deeper research for the future. In Section 3 we extend these results to the case of complex $\pp$, $|\pp-2|\leq1$.
The crucial consequence of these results is the fact that a function $\mathfrak{X}(\pp,x)$
(which gives a bijection between trees $\mathcal{Q}_{1}$ and
$\mathcal{Q}_{\pp}$) is a continuous function in $x$ and an analytic function
in $\pp$ for $|\pp-2|\leq 1$. In Section
4 we introduce exactly the same
integral transforms of $\fp(x)$ as was done in a special (though most
important) case of $F(x)=F_{1}(x)$. Also, in this section we prove certain
relations among the moments. In Section 5 we give the proof of the three term
functional equation for $G_{\pp}(z)$ and the integral equation for
$\mathfrak{m}_{\pp}(t)$. Finally, Theorem \ref{thm1.1} is proved in
Section 6. The hierarchy of sections is linear, and all results from previous
ones is used in Section 6. Appendix A. contains: derivation for the series (\ref{exr}); MAPLE codes to compute rational
functions $\mathbf{H}_{n}(z)$ and $\mathbf{Q}_{n}(z)$; description of high-precision method to calculate numerical
values for the constants $m_{L}$; auxiliary lemmas for the Section 3. The paper also contains
graphs of some $\pp-$Minkowski question mark functions $\fp(x)$ for real $\pp$,
and also pictures of locus points of elements of trees $\mathcal{Q}_{\pp}$ for
certain characteristic values of $\pp$.

\section{$\pp-$question mark functions and $\pp-$continued fractions}
In this section we introduce a family of natural generalizations of the Minkowski
question mark function $F(x)$. Let $1\leq \pp<2$. Consider the following binary
tree, which we denote by $\mathcal{Q}_{\pp}$. We start from the root $x=1$.
Further, each element (``root") $x$ of this tree generates two ``offsprings" by
the following rule:
\begin{eqnarray*}
x\mapsto \frac{\pp x}{x+1}, \quad \frac{x+1}{\pp}.
\end{eqnarray*}
We will use the notation $\mathcal{T}_{\pp}(x)=\frac{x+1}{\pp}$,
$\mathcal{U}_{\pp}(x)=\frac{\pp x}{x+1}$. When $\pp$ is fixed, we will
sometimes discard the subscript. Thus, the first four generations lead to
\begin{eqnarray}
\xymatrix @R=.5pc @C=.5pc { & & & & & {1\over 1} & & & & & & & \\
& & {\pp\over 2} \ar@{-}[urrr] & & & & & & {2\over \pp} \ar@{-}[ulll] & & & \\
& {\pp^{2}\over \pp+2} \ar@{-}[ur] & & {\pp+2\over 2\pp}\ar@{-}[ul] & & & & {2\pp\over \pp+2}\ar@{-}[ur] & & {\pp+2\over \pp^2} \ar@{-}[ul] & &\\
{\pp^3\over \pp^2+\pp+2} \ar@{-}[ur] & {\pp^2+\pp+2\over \pp^2+2\pp} \ar@{-}[u]
& {\pp^2+2\pp\over 3\pp+2} \ar@{-}[ur] & & {3\pp+2\over 2\pp^2} \ar@{-}[ul] & &
{2\pp^2\over 3\pp+2} \ar@{-}[ur] & & {3\pp+2\over \pp^2+2\pp} \ar@{-}[ul] &
{\pp^2+2\pp\over \pp^2+\pp+2} \ar@{-}[u] & {\pp^2+\pp+2\over \pp^3}
\ar@{-}[ul]}\label{cwp}
\end{eqnarray}
We refer the reader to the paper \cite{isola}, where authors consider a rather similar construction,
though having a different purpose in mind (see also \cite{bonnano}). Denote by $T_{n}(\pp)$ the sequence of polynomials,
appearing as numerators of
fractions of this tree. Thus, $T_{1}(\pp)=1$, $T_{2}(\pp)=\pp$, $T_{3}(\pp)=2$.
Directly from the definition of this tree we inherit that
\begin{eqnarray*}
T_{2n}(\pp)&=&\pp T_{n}(\pp) \text{ for }n\geq 1,\\
T_{2n-1}(\pp)&=&T_{n-1}(\pp)+\pp^{-\epsilon}T_{n}(\pp)\text{ for } n\geq 2,
\end{eqnarray*}
where $\epsilon=\epsilon(n)=1$ if $n$ is a power of two, and $\epsilon=0$ otherwise. Thus,
the definition of these polynomials is almost the same as it appeared in
\cite{kmp} (these polynomials were named Stern polynomials by the authors),
with the distinction that in \cite{kmp} everywhere one has $\epsilon=0$.
Naturally, this difference produces different sequence of polynomials.\\
\indent There are $2^{n-1}$ positive real numbers in each generation of the tree
$\mathcal{Q}_{\pp}$, say $a^{(n)}_{k}$, $1\leq k\leq 2^{n-1}$. Moreover, they
are all contained in the interval $[\pp-1,\frac{1}{\pp-1}]$. Indeed, this holds
for the initial root $x=1$, and
\begin{eqnarray*}
\pp-1\leq x\leq\frac{1}{\pp-1}\Leftrightarrow \pp-1\leq\frac{\pp x}{x+1}\leq
1,\\
\pp-1\leq x\leq\frac{1}{\pp-1}\Leftrightarrow 1\leq\frac{x+1}{\pp}\leq
\frac{1}{\pp-1}.
\end{eqnarray*}
This also shows that the left offspring is contained in the interval
$[\pp-1,1]$, while the right one - in the interval $[1,\frac{1}{\pp-1}]$. The
real numbers appearing in this tree have intrinsic relation with {\it
$\pp-$continued fractions algorithm.} The definition of the latter is as
follows. Let $x\in(\pp-1,\frac{1}{\pp-1})$. Consider the following procedure:
\begin{eqnarray*}
R_{\pp}(x)=\left\{\begin{array}{c@{\qquad}l} \mathcal{T}^{-1}(x)=\pp x-1, &
\mbox{if}\quad 1\leq x<\frac{1}{\pp-1},\\
\mathcal{I}(x)=\frac{1}{x}, & \mbox{if}\quad \pp-1< x< 1,\\
\text{STOP}, & \mbox{if}\quad x=\pp-1.
\end{array}\right.
\end{eqnarray*}
Then each such $x$ can be uniquely represented as $\pp-$continued fraction
\begin{eqnarray*}
x=[a_{0},a_{1},a_{2},a_{3},....]_{\pp},
\end{eqnarray*}
where $a_{i}\in\mathbb{N}$ for $i\geq 1$, and $a_{0}\in\mathbb{N}\cup\{0\}$.
This notation means that in the course of iterations $R_{\pp}^{\infty}(x)$ we
apply $\mathcal{T}^{-1}(x)$ exactly $a_{0}$ times, then once $\mathcal{I}$,
then we apply $\mathcal{T}^{-1}$ exactly $a_{1}$ times, then $\mathcal{I}$, and
so on. The procedure terminates exactly for those
$x\in(\pp-1,\frac{1}{\pp-1})$, which are the members of the tree
$\mathcal{Q}_{\pp}$ (``$\pp$-rationals"). Also, direct inspection shows that if
procedure does terminate, the last entry $a_{s}\geq 2$. Thus, we have the same
ambiguity for the last entry exactly as is the case with ordinary continued
fractions. At this point it is straightforward to show that the $n$th
generation of $\mathcal{Q}_{\pp}$ consists of $x=[a_{0},a_{1},...,a_{s}]_{\pp}$
such that
$\sum_{j=0}^{s}a_{j}=n$, exactly as in the case $\pp=1$ and tree (\ref{cw}).\\
\indent Now, consider a function $\mathfrak{X}_{\pp}(x)$ with the following property:
$\mathfrak{X}_{\pp}(x)=\overline{x}$, where $x$ is a rational number in the
Calkin-Wilf tree (\ref{cw}), and $\overline{x}$ is a corresponding number in
the tree (\ref{cwp}). In other words, $\mathfrak{X}_{\pp}(x)$ is simply a
bijection between these two trees. First, if $x<y$, then
$\overline{x}<\overline{y}$. Also, all positive rationals appear in the tree
(\ref{cw}) and they are everywhere dense in $\mathbb{R}_{+}$. Moreover,
$\mathcal{T}$ and $\mathcal{U}$ both preserve order, and
$[\pp-1,\frac{1}{\pp-1})$ is a disjoint union of
$\mathcal{T}[\pp-1,\frac{1}{\pp-1})$ and $\mathcal{U}[\pp-1,\frac{1}{\pp-1})$.
Now it is obvious that the function $\mathfrak{X}_{\pp}(x)$ can be extended to
a continuous monotone increasing function
\begin{eqnarray*}
\mathfrak{X}_{\pp}(\star):\mathbb{R}_{+}\rightarrow
[\pp-1,\frac{1}{\pp-1}\Big{)}, \quad\mathfrak{X}_{\pp}(\infty)=\frac{1}{\pp-1}.
\end{eqnarray*}
Thus,
\begin{eqnarray*}
\mathfrak{X}_{\pp}\big{(}[a_{0},a_{1},a_{2},a_{3}...]\big{)}=[a_{0},a_{1},a_{2},a_{3}...]_{\pp}.
\end{eqnarray*}
As can be seen from the definitions of both trees (\ref{cw}) and (\ref{cwp}),
this function satisfies functional equations
\begin{eqnarray}
\mathfrak{X}_{\pp}(x+1)&=&\frac{\mathfrak{X}_{\pp}(x)+1}{\pp},\nonumber\\
\mathfrak{X}_{\pp}\Big{(}\frac{x}{x+1}\Big{)}&=&\frac{\pp\mathfrak{X}_{\pp}(x)}{\mathfrak{X}_{\pp}(x)+1},\label{sing}\\
\mathfrak{X}_{\pp}\Big{(}\frac{1}{x}\Big{)}&=&\frac{1}{\mathfrak{X}_{\pp}(x)}.\nonumber
\end{eqnarray}
The last one (symmetry property) is a consequence of the first two. We are not
aware whether this notion of $\pp-$continued fractions is new or not. For
example,
\begin{eqnarray*}
\frac{1+\sqrt{1+4\pp}}{2\pp}&=&[1,1,1,1,1,1,...]_{\pp}=\mathfrak{X}_{\pp}\Big{(}\frac{1+\sqrt{5}}{2}\Big{)},\\
\sqrt{3}&=&[4,2,1,10,1,1,2,1,5,1,1,2,1,2,1,1,2,1,3,7,4,...]_{\frac{3}{2}},\\
2&=&[4,1,1,\overline{2,1,1}]_{\sqrt{2}}.
\end{eqnarray*}
Now fix $\pp$, $1\leq \pp<2$. The following proposition follows immediately
from the properties of $F(x)$.
\begin{prop}There exists a limit distribution of the $n$th generation of the tree $\mathcal{Q}_{\pp}$ as
$n\rightarrow\infty$, defined as
\begin{eqnarray*}
\fp(x)=\lim\limits_{n\rightarrow\infty}2^{-n+1}\#\{k:a^{(n)}_{k}<x\}.
\end{eqnarray*}
This function is continuous, $\fp(x)=0$ for $x\leq \pp-1$, $\fp(x)=1$ for
$x\geq \frac{1}{\pp-1}$, and it satisfies two functional equations:
\begin{eqnarray}
2\fp(x)=\left\{\begin{array}{c@{\qquad}l} \fp(\pp x-1)+1, & \mbox{if}\quad
1\leq x\leq \frac{1}{\pp-1},
\\ \fp({x\over \pp-x}), & \mbox{if}\quad \pp-1\leq x\leq 1. \end{array}\right.
\label{distrp}
\end{eqnarray}
Additionally,
\begin{eqnarray*}
\fp(x)+\fp\Big{(}\frac{1}{x}\Big{)}=1\text{ for }x>0.
\end{eqnarray*}
The explicit expression for $\fp(x)$ is given by
\begin{eqnarray*}
\fp([a_{0},a_{1},a_{2},a_{3},...]_{\pp})=1-2^{-a_{0}}+2^{-(a_{0}+a_{1})}-2^{-(a_{0}+a_{1}+a_{2})}+....
\label{minp}
\end{eqnarray*}
\label{prop2.1}
\end{prop}
We will refer to the last functional equation as {\it the symmetry property}.
As was said, it is a consequence of the other two, though it is convenient to
separate it.\\

{\it Proof. } Indeed, as it is obvious from the observations above, we simply
have
\begin{eqnarray*}
\fp\big{(}\mathfrak{X}_{\pp}(x)\big{)}=F(x),\quad x\in[0,\infty).
\end{eqnarray*}
Therefore, two functional equations follow from (\ref{distr}) and (\ref{sing}).
All the other statements are immediate and follow from the
properties of $F(x)$. $\square$\\

Equally important, consider the binary tree (\ref{cwp}) for
$\pp>2$. In this case analogous proposition holds.
\begin{prop}Let $\pp>2$. Then there exists a limit distribution of the $n$th generation as
$n\rightarrow\infty$. Denote it by $f_{\pp}(x)$ This function is continuous,
$f_{\pp}(x)=0$ for $x\leq \frac{1}{\pp-1}$, $f_{\pp}(x)=1$ for $x\geq \pp-1$,
and it satisfies two functional equations:
\begin{eqnarray*}
2f_{\pp}(x)=\left\{\begin{array}{c@{\qquad}l} f_{\pp}(\pp x-1) & \mbox{if}\quad 1\leq
x\leq \pp-1,
\\ f_{\pp}({x\over \pp-x})+1 & \mbox{if}\quad \frac{1}{\pp-1}\leq x\leq 1, \end{array}\right.
\end{eqnarray*}
and
\begin{eqnarray*}
f_{\pp}(x)+f_{\pp}\Big{(}\frac{1}{x}\Big{)}=1\text{ for }x>0.
\end{eqnarray*}
\label{prop2.2}
\end{prop}
{\it Proof. }The proof is analogous to the one of Proposition \ref{prop2.1}, only this time
we use equivalences
\begin{eqnarray*}
\pp-1\leq x\leq\frac{1}{\pp-1}\Leftrightarrow 1\leq\frac{\pp x}{x+1}\leq
\pp-1,\\
\pp-1\leq x\leq\frac{1}{\pp-1}\Leftrightarrow
\frac{1}{\pp-1}\leq\frac{x+1}{\pp}\leq \pp-1.\quad\square
\end{eqnarray*}
For the sake of uniformity, we introduce $\fp(x)=1-f_{\pp}(x)$ for $\pp>2$.
Then $\fp(x)$ satisfies exactly the same functional equations (\ref{distr}),
with a slight difference that $\fp(x)=1$ for $x\leq\frac{1}{\pp-1}$ and
$\fp(x)=0$ for $x\geq\pp-1$. Consequently, we will not separate these two cases
and all our subsequent results hold uniformly. To this account it should be
noted
 that, for example, in case $\pp>2$ the integral $\int_{\pp-1}^{1}\star\d\star$ should be understood as
 $-\int_{1}^{\pp-1}\star\d\star$. Figure 1 gives graphic images of typical cases for $\fp(x)$.\\
\begin{figure}
\centering
\begin{tabular}{c c}

\epsfig{file=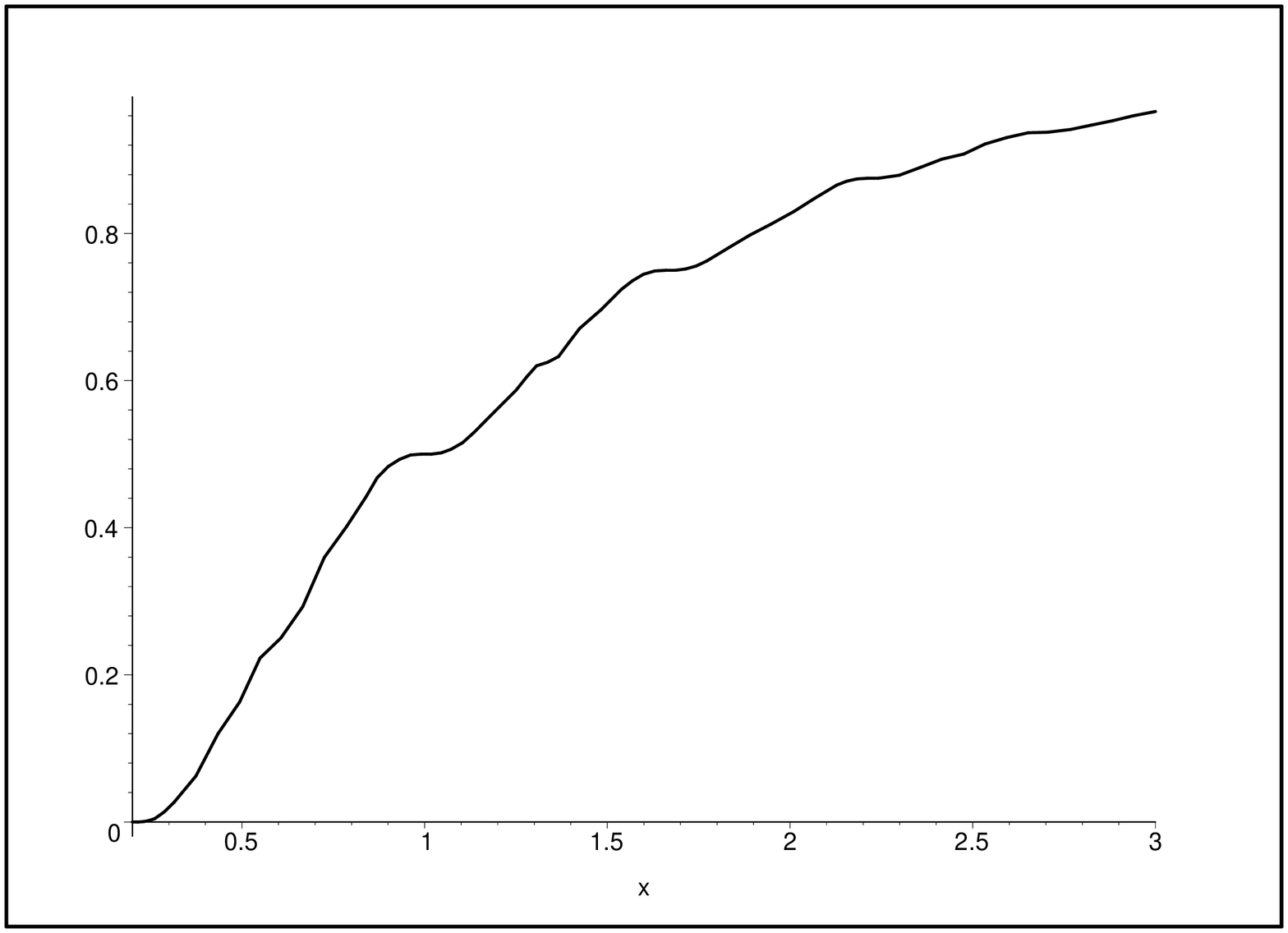,width=188pt,height=230pt,angle=-90}
&\epsfig{file=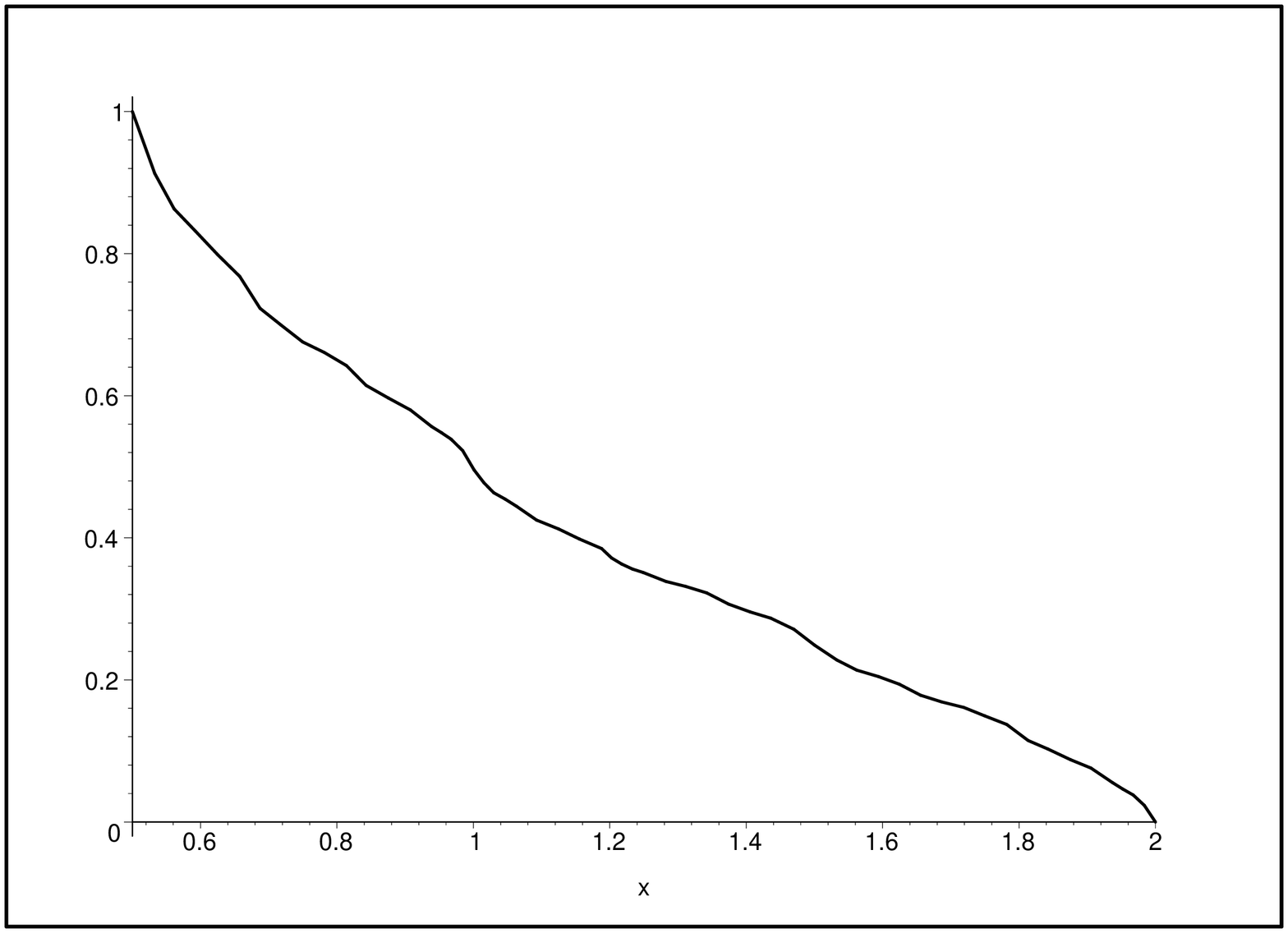,width=188pt,height=230pt,angle=-90} \\
$\pp=1.2$, $x\in[0.2,3]$ & $\pp=3$, $x\in[0.5,2]$\\

\epsfig{file=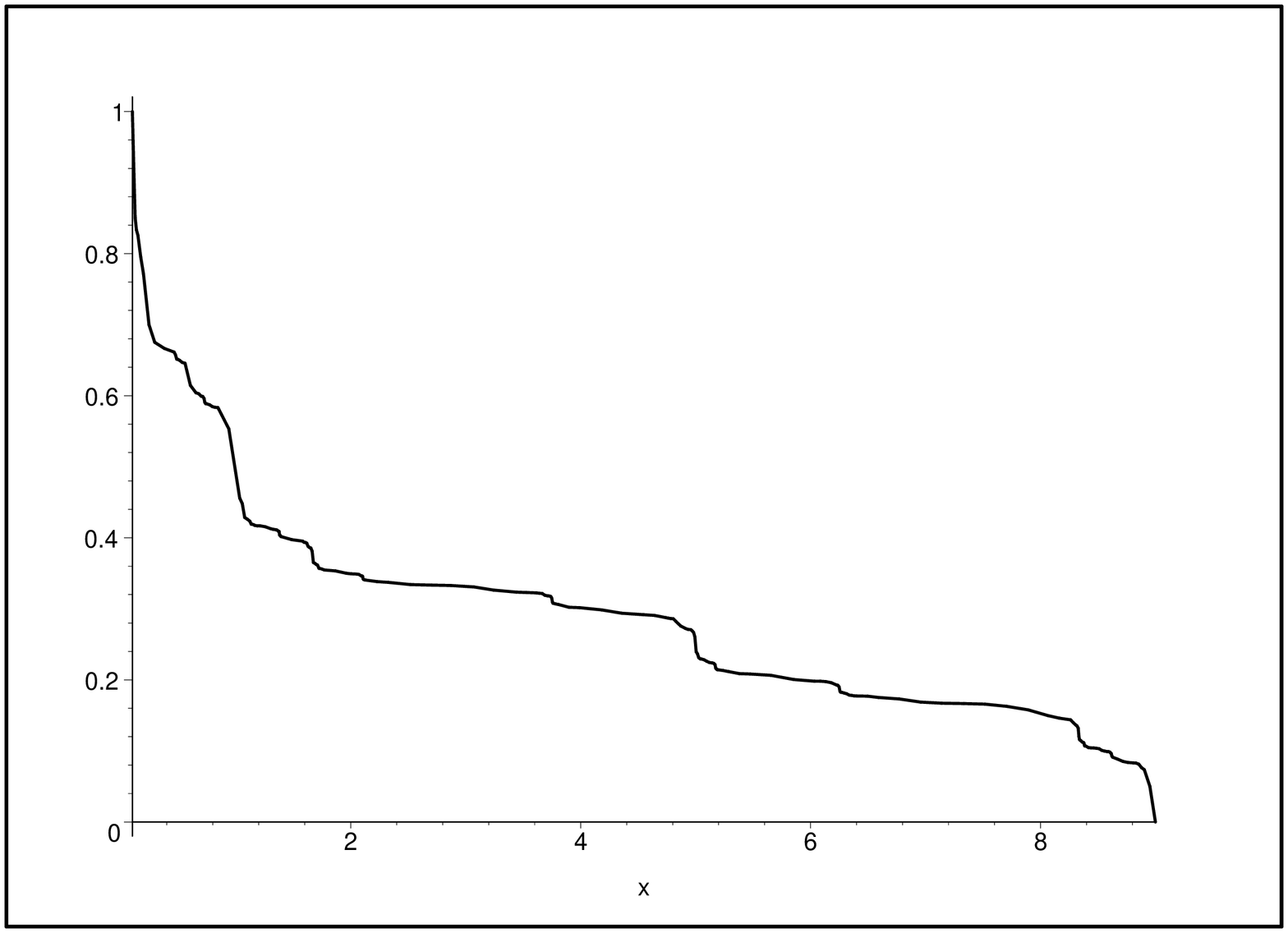,width=188pt,height=230pt,angle=-90} &
\epsfig{file=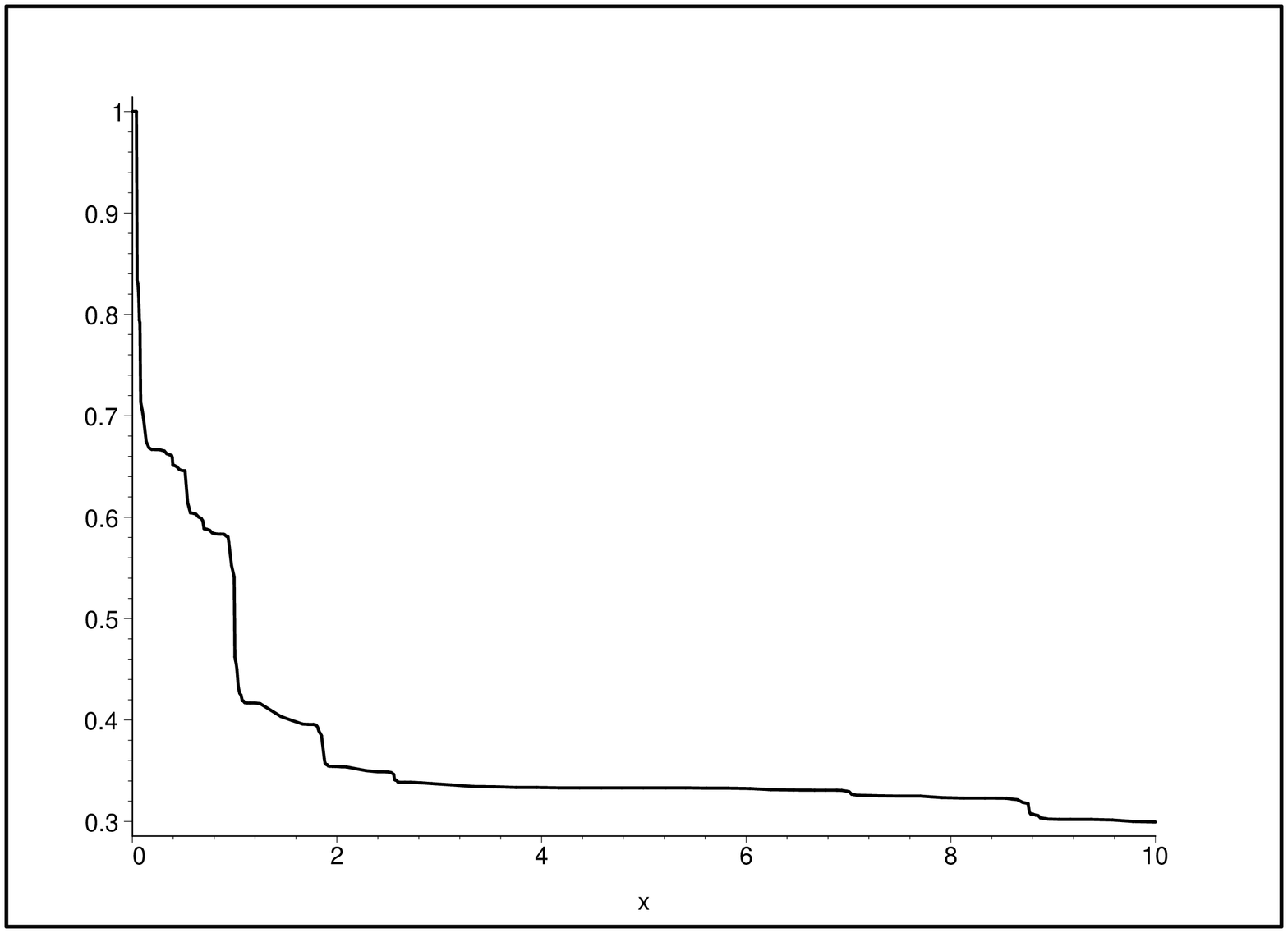,width=188pt,height=230pt,angle=-90}\\
$\pp=10$, $x\in[0.1,9]$ & $\pp=25$, $x\in[0,10]$\\
\end{tabular}
\caption{Functions $\fp(x)$}.
\end{figure}
\section{Complex case}
After dealing the case of real $\pp$, $1\leq \pp<\infty$, let us consider a
tree (\ref{cwp}) when $\pp\in\mathbb{C}$. For our purpose we will concentrate on the case $|\pp-2|\leq 1$. It should be noted that the method
which we use allows to extend these result to the case $\Re\pp\geq 1$. The question in determining the set in the complex plain where
similar results are valid remains open. More importantly, the problem to determine all $\pp\in\mathbb{C}$ for which there exists an analytic
function $G_{\pp}(z)$, which satisfied (\ref{fep}), seems to be much harder and interesting. Here and below
$[0,\infty]$ stands for a compactification of $[0,\infty)$. In the sequel, the notion of a function $f(z)$ to be analytic in the closed disc
$|z-2|\leq 1$ means that for $z_{0}\neq 1$, $|z_{0}-2|\leq 1$, this function is analytic in a certain small neighborhood of $z_{0}$.
If $z_{0}=1$, this means that there exist all higher derivatives, if one approaches the point $z_{0}=1$ while remaining in the disc $|z-2|\leq 1$.\\
\noindent In this section we prove the following result.
\begin{thm}
There exists a unique function $\mathfrak{X}_{\pp}(x)=\mathfrak{X}(\pp,x):\{|\pp-2|\leq 1\}\times \{[0,\infty]\}\rightarrow\mathbb{C}\cup\{\infty\}$,
having these properties:\\
$(i)$ $\mathfrak{X}(\pp,x)$ satisfies functional equations (\ref{sing});\\
$(ii)$ For fixed $\pp\neq 1$, $\mathfrak{X}(\pp,x):[0,\infty]\rightarrow\mathbb{C}$ is a continuous function, and the image (denote it by $\I_{\pp}$) is thus a bounded curve;
it is contained in the domain $\{\mathbb{C}\setminus\{|z+1|\leq \frac{3}{4}\}$;\\
$(iii)$ For every $\pp$, $|\pp-2|\leq 1$, $\pp\neq 1$, in some neighborhood of $\pp$ there exists the derivative
$\frac{\p}{\p\pp}\mathfrak{X}(\pp,x)$, which is a continuous function for $x\in[0,\infty]$;\\
$(iv)$ There exist all derivatives
$\mathcal{S}_{N}(x)=\frac{\p^{N}}{\p\pp^{N}}\mathfrak{X}(\pp,x)|_{\pp=1}:[0,\infty)\rightarrow\mathbb{R}$
(the derivatives are taken inside $|\pp-2|\leq 1$). These
functions are uniformly continuous for irrational $x$ in any finite interval.
Moreover, $\mathcal{S}_{N}(x)\ll_{N}x^{N+1}$ for $x\geq 1$, and $\mathcal{S}_{N}(x)\ll_{N}1$ for $x\in(0,1)$.
\label{thm3.1}
\end{thm}
The curve $\I_{\pp}$ has a natural fractal structure: it decomposes into two
parts, namely $\frac{\I_{\pp}+1}{\pp}$ and $\frac{\pp\I_{\pp}}{\I_{\pp}+1}$,
with a common point $z=1$. Additionally, $\I_{\pp}=\frac{1}{\I_{\pp}}$. As a consequence, $0\notin\I_{\pp}$ for $\pp\neq 1$.
Figures 2-4 show the images of $\I_{\pp}$ for certain characteristic values of $\pp$.\\

The investigations of the tree $\mathcal{Q}_{\pp}$ deserve a separate paper. I am very grateful to my colleagues Jeffrey Lagarias
and Stefano Isola, who sent me various references, also informing about the intrinsic relations of this problem with: Julia sets of
rational maps of the Riemann sphere; iterated function systems; forward limit sets of semigroups; various topics from complex dynamics
and geometry of discrete groups. Thus, the problem is much more subtle and involved than it appears to be. This poses a difficult question
on the limit set of the semigroup generated by transformations $\mathcal{U}_{\pp}$ and $\mathcal{T}_{\pp}$, or any other two ``conjugate"
analytic maps of the Riemann sphere (say, two analytic maps $\mathcal{A}$ and $\mathcal{B}$ are ``conjugate", if $\mathcal{A}(\alpha)=\alpha$,
$\mathcal{B}(\beta)=\beta$, $\mathcal{A}(\beta)=\mathcal{B}(\alpha)$ for some two points $\alpha$ and $\beta$ on the Riemann
sphere). Possibly, certain techniques from complex dynamics do apply here. As pointed out by Curtis McMullen, the property of boundedness of
$\I_{\pp}$ can be reformulated in a coordinate-free manner. It appears that this curve consists of the closure of the attracting fixed points
of the elements of the semigroup $\langle\mathcal{T}_{\pp},\mathcal{U}_{\pp}\rangle$. Then the property for the curve being
bounded and being bounded away from $z=0$ means that it does not contain a repelling fixed point of $\mathcal{T}_{\pp}$ ($z=\infty$)
and a repelling fixed point
of $\mathcal{U}_{\pp}$ ($z=0$). It also does not contain neither of the repelling fixed points of the elements of this semigroup. Note that
$\mathcal{T}_{2}(1)=\mathcal{U}_{2}(1)=1$, $\mathcal{T}'_{2}(1)=\mathcal{U}'_{2}(1)=1/2$. Thus, there exists a small ball $\mathbf{D}$ around $z=1$, such
that $\mathcal{T}_{2}(\mathbf{D})\subset\mathbf{D}$,  $\mathcal{U}_{2}(\mathbf{D})\subset\mathbf{D}$, and the last two maps are contractions in
$\mathbf{D}$. This strict containment is an open condition on $\pp$, and thus there exists a neighborhood of $\pp=2$ such that Theorem \ref{thm3.1}
does hold. I am grateful to Curtis McMullen for this remark: we get the result almost for free. Yet, the full result for $|\pp-2|\leq 1$ is needed.
This is not a new kind of problem. Some cases of pairs of M\"{o}bius transformations
were studied. For example, the author in \cite{bousch} deals with the case of a semigroup generated by two maps
$z\mapsto sz\pm 1$, for fixed $s$, $|s|<1$, and investigates a closure of a set of all attracting fixed points. For example, for
$|s|>2^{-1/2}$ this set is connected. Further development of this problem
can be seen in \cite{solomyak}. On the other hand, the case of one rational map is rather well understood, and it is treated in \cite{beardon}.
Thus, though the machinery of complex dynamics can greatly clarify our understanding of the structure of the curve $\I_{\pp}$, we will rather employ the techniques
 from the analytic theory of continued fractions. The main source is the monograph by H.S. Wall \cite{wall}. (Lemmas \ref{lema.1},
 \ref{lema.2} and \ref{lema.3} can be found in the Appendix A.2.)\\

\noindent {\bf Proof of Theorem \ref{thm3.1}}. We need the following two results.
\begin{thm}(\cite{wall} p. 57.) Let $v_{\nu}$, $\nu\in\mathbb{N}$ be positive numbers such that
\begin{eqnarray}
v_{1}<1,\quad v_{\nu}+v_{\nu+1}\leq 1,\text{ for }\nu\geq 1.
\label{sal1}
\end{eqnarray}
Suppose given complex numbers $e_{\nu}$, $\nu\in\mathbb{N}$, such that
\begin{eqnarray}
|e_{\nu+1}|-\Re(e_{\nu+1})\leq v_{\nu},\quad\nu\geq 1.
\label{sal2}
\end{eqnarray}
Define the sequence $b_{\nu}$ by the recurrence $b_{1}=1$, $e_{\nu+1}=\frac{1}{b_{\nu}b_{\nu+1}}$, $\nu\geq 1$. Then the continued fraction
\begin{eqnarray}
\mathcal{F}=\cfrac{1}{1+\cfrac{e_{2}}{1+\cfrac{e_{3}}{1+\cfrac{e_{4}}{\ddots}}}}
\label{bas}
\end{eqnarray}
converges if, and only if, $(a)$ some $e_{\nu}$ vanishes, or $(b)$ $e_{\nu}\neq 0$ for $\nu\geq 2$ and the series $\sum_{\nu=1}^{\infty}|b_{\nu}|$
 diverges.
Moreover, if $e_{\nu}(z):\mathbf{K}_{1}\rightarrow\mathbf{K}_{2}$ are analytic functions of a complex variable, $\mathbf{K}_{1}$ and $\mathbf{K}_{2}$ are
compact sets, (\ref{sal1}) and (\ref{sal2}) are satisfied, and the above series diverges uniformly, then the continued fraction converges uniformly for all $z\in\mathbf{K}_{1}$.
\label{thm3.2}
\end{thm}
\begin{thm}(\cite {wall}, p. 60.) If all $v_{\nu}=\frac{1}{2}$, and the conditions $(a)$ and $(b)$ of Theorem \ref{thm3.2} hold,
then $|\mathcal{F}-1|\leq 1$, $\mathcal{F}\neq 0$.
\label{thm3.3}
\end{thm}
For $a,b\in\mathbb{N}$, $\pp\in\mathbb{C}$, $|\pp-2|\leq 1$, define rational functions
\begin{eqnarray*}
W_{a}(\pp)&=&\frac{\pp^{a}-1}{\pp^{a+1}-\pp^{a}},\\
T_{a,b}(\pp)&=&W^{-1}_{a}(\pp)W^{-1}_{b}(\pp)\pp^{-a}=\frac{(\pp-1)^{2}\pp^{b}}{(\pp^{a}-1)(\pp^{b}-1)},\quad
T_{a,\infty}(\pp)=\frac{(\pp-1)^{2}}{(\pp^{a}-1)}.
\end{eqnarray*}
Since, for fixed $\pp\neq 1$, $W_{a}(\pp)\rightarrow\pp-1$, as $a\rightarrow\infty$,
then there exist two constants $k_{1}=k_{1}(\pp)$ and $k_{2}=k_{2}(\pp)$, such that
\begin{eqnarray}
0<k_{1}\leq|W_{a}(\pp)|\leq k_{2}<+\infty,\quad a\in\mathbb{N}.\label{wpol}
\end{eqnarray}
Let $x\geq 1$, $x=[a_{1},a_{2},a_{3},...]$, be an irrational number, $a_{i}\in\mathbb{N}$. Let us consider the continued fraction
\begin{eqnarray}
\mathcal{F}(\pp,x)=\mathcal{F}(\pp,a_{1},a_{2},...)=
\cfrac{1}{1+\cfrac{T_{a_{1},a_{2}}(\pp)}{1+\cfrac{T_{a_{2},a_{3}}(\pp)}{1+\cfrac{T_{a_{3},a_{4}}(\pp)}{\ddots}}}}.
\label{cfrac}
\end{eqnarray}
If $x=[a_{1},a_{2},...,a_{\kappa}]\geq 1$ is rational, let us define
\begin{eqnarray*}
\mathcal{F}(\pp,x)=\mathcal{F}(\pp,a_{1},a_{2},...,a_{\kappa})=
\cfrac{1}{1+\cfrac{T_{a_{1},a_{2}}(\pp)}{1+\cfrac{T_{a_{2},a_{3}}(\pp)}{1+\cfrac{\ddots}{1+T_{a_{\kappa},\infty}}}}}.
\end{eqnarray*}
From the definition, this continued fraction obeys the following rule
\begin{eqnarray*}
\mathcal{F}(\pp,a_{1},a_{2},...)=\frac{1}{1+T_{a_{1},a_{2}}(\pp)\cdot\mathcal{F}(\pp,a_{2},a_{3}...)}.
\end{eqnarray*}
We will now apply Theorem \ref{thm3.2} to $\mathcal{F}(\pp,a_{1},a_{2},a_{3},...)$. Suppose $x$ is irrational.
Thus, let $e_{\nu}=T_{a_{\nu-1},a_{\nu}}(\pp)$, $\nu\geq 2$.
Let us define constants
\begin{eqnarray*}
\mu(a,b)=\sup\limits_{\pp\in\mathbb{C},|\pp-2|\leq 1}|T_{a,b}(\pp)|-\Re(T_{a,b}(\pp)).
\end{eqnarray*}
By Lemma \ref{lema.1}, $\mu(a,b)+\mu(b,c)<0.76$, $a,b,c\in\mathbb{N}$. Further,
from the definition in Theorem \ref{thm3.2} it follows that
\begin{eqnarray}
b_{2\nu}=W_{a_{1}}(\pp)W_{a_{2\nu}}(\pp)\pp^{a_{2\nu-1}-...+a_{3}-a_{2}+a_{1}},\nonumber\\
b_{2\nu+1}=W^{-1}_{a_{1}}(\pp)W_{a_{2\nu+1}}(\pp)\pp^{a_{2\nu}-...-a_{3}+a_{2}-a_{1}}.\label{bb}
\end{eqnarray}
It is obvious that the series $\sum_{\nu=1}^{\infty}|b_{\nu}|$ diverges.
Hence, Theorem \ref{thm3.2} tells that the continued fraction converges, and that for fixed irrational $x=[a_{1},a_{2},...]>1$,
$\mathcal{F}(\pp_{0},a_{1},a_{2},...)$ is an analytic function in $\pp_{0}$ in some small neighborhood of $\pp$.
For rational $x$ this is in fact a rational function.\\
\indent As it is shown in \cite{wall}, the $\nu$th convergent of the continued
fraction (\ref{bas}) (denote it by $\frac{A_{\nu}}{B_{\nu}}$) is equal to the $\nu$th convergent (denote it by $\frac{P_{\nu}}{Q_{\nu}}$) of the continued fraction
\begin{eqnarray*}
\cfrac{1}{b_{1}+\cfrac{1}{b_{2}+\cfrac{1}{b_{3}+\cfrac{1}{\ddots}}}}.
\end{eqnarray*}
Moreover, since (\ref{sal1}) and (\ref{sal2}) are satisfied, we have that, for certain positive constant $k=k(b_{1},b_{2},b_{3})$ (\cite{wall}, p.55-56),
\begin{eqnarray}
|Q_{2\nu}|&\geq& k(1+|b_{2}|+|b_{4}|+...,+|b_{2\nu}|),\nonumber\\
|Q_{2\nu+1}|&\geq& k(1+|b_{3}|+|b_{5}|+...,+|b_{2\nu+1}|),\label{ineq}\\
\Big{|}\frac{A_{\nu+1}}{B_{\nu+1}}-\frac{A_{\nu}}{B_{\nu}}\Big{|}&=&\frac{1}{|Q_{\nu}Q_{\nu+1}|}.\nonumber
\end{eqnarray}
Now we have
\begin{prop} Fix $\pp\in\mathbb{C}$, $|\pp-2|\leq 1$, $\pp\neq 1$.
Let $x=[a_{1},a_{2},...]\geq 1$ be a real number. The function $\mathcal{F}(\pp,x):[1,\infty]\rightarrow\mathbb{C}$ is continuous.
\label{prop3.4}
\end{prop}
\proof
Fix irrational $x>1$. Let $\delta>0$, and $y\geq 1$ be such that $|x-y|<\delta$. Then there exists $N$ such
that the first $N$ partial quotients of $x$ and $y$ coincide,
$N=N(\delta)\rightarrow\infty$ as $\delta\rightarrow 0$. Consequently,
let the corresponding convergents to $\mathcal{F}(\pp,x)$ and $\mathcal{F}(\pp,y)$ be respectively
\begin{eqnarray*}
&&\frac{A_{1}}{B_{1}},\quad \frac{A_{2}}{B_{2}},\quad ...,\frac{A_{N}}{B_{N}},\quad\frac{A_{N+1}}{B_{N+1}},\quad\frac{A_{N+2}}{B_{N+2}},...;\text{ and}\\
&&\frac{A_{1}}{B_{1}},\quad \frac{A_{2}}{B_{2}},\quad ...,\frac{A_{N}}{B_{N}},\quad\frac{A'_{N+1}}{B'_{N+1}},\quad\frac{A'_{N+2}}{B_{'N+2}}...
\end{eqnarray*}
Now, combining (\ref{wpol}), (\ref{bb}) and (\ref{ineq}) we see that
\begin{eqnarray*}
|Q_{2\nu}Q_{2\nu+1}|> k^{2}k^{3}_{1}k^{-1}_{2}&\times&\Big{(}|\pp|^{a_{1}}+|\pp|^{a_{3}-a_{2}+a_{1}}+...+|\pp|^{a_{2\nu-1}-...+a_{3}-a_{2}+a_{1}}\Big{)}\\
&\times&\Big{(}|\pp|^{a_{2}-a_{1}}+|\pp|^{a_{4}-a_{3}+a_{2}-a_{1}}+...+|\pp|^{a_{2\nu}-...-a_{3}+a_{2}-a_{1}}\Big{)}.
\end{eqnarray*}
Denote $c_{1}=k^{2}k^{3}_{1}k^{-1}_{2}$. Let $|\pp|^{a_{2\ell-1}-...+a_{3}-a_{2}+a_{1}}=\lambda_{\ell}$, $1\leq \ell\leq \nu$. The above inequality and
the arithmetic-harmonic mean inequality give
\begin{eqnarray*}
|Q_{2\nu}Q_{2\nu+1}|&>& c_{1}(\lambda_{1}+\lambda_{2}+...+\lambda_{\nu})\cdot
(|\pp|^{a_{2}}\lambda^{-1}_{1}+|\pp|^{a_{4}}\lambda_{2}^{-1}+...+|\pp|^{a_{2\nu}}\lambda^{-1}_{\nu})\\
&\geq&|\pp|c_{1}(\lambda_{1}+\lambda_{2}+...+\lambda_{\nu})\cdot(\lambda^{-1}_{1}+\lambda^{-1}_{2}+...+\lambda^{-1}_{\nu})
\geq |\pp|c_{1}\nu^{2}\quad \nu\geq 1.
\end{eqnarray*}
Analogously we prove that $|Q_{2\nu-1}Q_{2\nu}|>|\pp|c_{2}\nu^{2}$, $\nu\geq 2$. Thus, $|Q_{\nu}Q_{\nu+1}|>c\nu^{2}$ for certain
real $c>0$, $\nu\geq 2$. We see that (\ref{ineq}) yield
\begin{eqnarray*}
\Big{|}\mathcal{F}(\pp,x)-\frac{A_{N}}{B_{N}}\Big{|}<\sum\limits_{\nu=N}^{\infty}\frac{1}{|Q_{\nu}Q_{\nu+1}|}\leq
\sum\limits_{\nu=N}^{\infty}\frac{c^{-1}}{\nu^{2}}<\frac{c^{-1}}{N-1}; \quad \Big{|}\mathcal{F}(\pp,y)-\frac{A_{N}}{B_{N}}\Big{|}<\frac{c^{-1}}{N-1}.
\end{eqnarray*}
This implies $|\mathcal{F}(\pp,x)-\mathcal{F}(\pp,y)|<\frac{2c^{-1}}{N-1}$. In case  $x$ is rational we argue in a similar way.
In this case note that real numbers close to $x=[a_{1},a_{2},...,a_{\kappa}]$ are of the form or $[a_{1},a_{2},...,a_{\kappa},T,...]$, either
$[a_{1},a_{2},...,a_{\kappa}-1,1,T,...]$ for $T$ sufficiently large. The case $x=\infty$ is analogous. This establishes the validity of the Proposition. $\square$\\

Eventually, for real number $x\geq 0$, $x=[a_{0},a_{1},a_{2},...]$, let us define
\begin{eqnarray*}
\mathfrak{X}(\pp,[a_{0},a_{1},...])=W_{a_{0}}(\pp)+\cfrac{\pp^{-a_{0}}}{W_{a_{1}}(\pp)+\cfrac{\pp^{-a_{1}}}{W_{a_{2}}(\pp)+\cfrac{\pp^{-a_{2}}}
{W_{a_{3}}(\pp)+\ddots}}}.
\end{eqnarray*}
After an equivalence transformation (\cite{wall}, p.19), this can be given an expression
\begin{eqnarray}
\mathfrak{X}(\pp,[a_{0},a_{1},...])=W_{a_{0}}(\pp)+\pp^{-a_{0}}W_{a_{1}}^{-1}(\pp)\cdot\mathcal{F}(\pp,a_{1},a_{2},a_{3},...).
\label{event}
\end{eqnarray}
From the very construction, this function satisfies the functional equations (\ref{sing}), is continuous at $x=1$ and thus is continuous
for $x\in[0,\infty]$. Obviously, (\ref{sing}) determine the values of $\mathfrak{X}(\pp,x)$ at rational $x$ uniquely,
hence a continuous solution to (\ref{sing}) is unique. We are left to show that the image of the curve
$\I_{\pp}$ is contained outside the circle $|z+1|\leq \frac{3}{4}$. This is equivalent to the statement that $\frac{\pp\I_{\pp}}{\I_{\pp}+1}$ is contained
inside the circle $|z-\pp|\leq\frac{4\pp}{3}$. But the points on $\frac{\pp\I_{\pp}}{\I_{\pp}+1}$ are exactly the point on the curve $\I_{\pp}$ with $a_{0}=0$.
Thus, we need to show that
\begin{eqnarray}
|\pp^{-1}\mathfrak{X}(\pp,[0,a_{1},a_{2},...])-1|=|\pp^{-1}W^{-1}_{a_{1}}\mathcal{F}(\pp,a_{1},a_{2},...)-1|\leq \frac{4}{3}.
\label{srytis}
\end{eqnarray}
Unfortunately, we cannot apply Theorem \ref{thm3.3} directly to all $\pp$, $|\pp-2|\leq 1$, since the table above Lemma \ref{lema.1}
shows that $\mu(1,b)>\frac{1}{2}$ for infinitely many $b$. The maximum values $\mu(1,b)$ (see the definition of this constant) are produced by points
$\pp$ close to $\chi=2+e^{2\pi i/3}$, or to $\overline{\chi}$. For this reason let us introduce
\begin{eqnarray*}
\mu^{\star}(a,b)=\sup\limits_{\pp\in\mathbb{C},|\pp-2|\leq 1,|\pp-\chi|\geq 0.19,|\pp-\overline{\chi}|\geq 0.19}|T_{a,b}(\pp)|-\Re(T_{a,b}(\pp)).
\end{eqnarray*}
Then indeed $\mu^{\star}(a,b)<\frac{1}{2}$ for all $a,b\in\mathbb{N}$. Thus, Theorem \ref{thm3.3} gives
$|\mathcal{F}(\pp,a_{1},a_{2},...)-1|\leq 1$, and the statement (\ref{srytis}) follows from Lemma \ref{lema.3}. In case $|\pp-2|\leq 1$,
$|\pp-\chi|<0.19$ (or $|\pp-\overline{\chi}|<0.19$) we use another theorem by Wall (\cite{wall}, p. 152), which describes the value region of a
continued fraction (\ref{bas}), provided elements $e_{\nu}$
belong to the compact domain in the parabolic region $|z|-\Re(ze^{i\phi})\leq 2h\cos^{2}\frac{\phi}{2}$,
for certain fixed $-\pi<\phi<+\pi$, $0<h\leq \frac{1}{4}$. We omit the details.
This proves part $(ii)$. In a similar fashion we prove part $(iii)$.
Finally, a direct inspection shows that slightly modified proofs remain valid in case $\pp=1$, if we define a function to be analytic at $\pp=1$, if it
possesses all higher $\pp-$derivatives, while remaining inside the disc $|\pp-2|\leq 1$. $\square$
\begin{defin} We define Minkowski $\pp-$question mark function
$\fp(x):\I_{\pp}\rightarrow[0,1]$, by
\begin{eqnarray*}
\fp(\mathfrak{X}(\pp,x))=F(x),\quad x\in[0,\infty].
\end{eqnarray*}
\label{defi3.5}
\end{defin}
\section{Properties of integral transforms of $\fp(x)$}
For given $\pp$, $|\pp-2|\leq 1$, we define
\begin{eqnarray*}
\chi_{n}=\frac{\pp+\pp^{n-1}-2}{\pp^{n-1}(\pp-1)},\quad
\I_{n}=[\chi_{n},\chi_{n+1}]=\mathfrak{X}(\pp,[n,n+1])\text{ for
}n\in\mathbb{N}_{0}.
\end{eqnarray*}
Complex numbers $\chi_{n}$ stand for the analogue of non-negative integers on
the curve $\I_{\pp}$. In other words, $\chi_{n}=\mathcal{U}^{n}(\pp-1)$. We
consider $\I_{n}$ as part of the curve $\I_{\pp}$ contained between the points
$\chi_{n}$ and $\chi_{n+1}$. Thus, $\chi_{0}=\pp-1$, $\chi_{1}=1$, and the
sequence $\chi_{n}$ is ``increasing", in the sense that $\chi_{j}$ as a point
on a curve $\I_{\pp}$ is between $\chi_{i}$ and $\chi_{k}$ if $i<j<k$.
Moreover,
$\bigcup\limits_{n=0}^{\infty}\I_{n}\bigcup\{\frac{1}{\pp-1}\}=\I_{\pp}$.
\begin{prop}
Let $\omega(x):\I_{\pp}\rightarrow\mathbb{C}$ be a continuous function. Then
\begin{eqnarray*}
\int\limits_{\I_{\pp}}\omega(x)\d\fp(x)=\sum\limits_{n=0}^{\infty}\frac{1}{2^{n+1}}\int\limits_{\I_{\pp}}
\omega\Big{(}\frac{x}{\pp^{n-1}(x+1)}+\frac{\pp^{n}-1}{\pp^{n+1}-\pp^{n}}\Big{)}\d\fp(x).
\end{eqnarray*}
\label{prop4.1}
\end{prop}
{\it Proof.} Indeed, using (\ref{distrp}) we obtain
\begin{eqnarray*}
\int\limits_{\I_{\pp}}\omega(x)\d\fp(x)=\sum\limits_{n=0}^{\infty}\int\limits_{\I_{n}}\omega(x)\d\fp(x)=
\sum\limits_{n=0}^{\infty}\int\limits_{\mathcal{T}^{n}(\I_{0})}\omega(x)\d\fp(x)\mathop{=}^{x\rightarrow \mathcal{T}^{n}x}\\
\sum\limits_{n=0}^{\infty}\frac{1}{2^{n}}\int\limits_{\I_{0}}\omega(\mathcal{T}^{n}x)\d\fp(x)\mathop{=}^{x\rightarrow
\mathcal{U}x}
\sum\limits_{n=0}^{\infty}\frac{1}{2^{n+1}}\int\limits_{\I_{\pp}}\omega(\mathcal{T}^{n}\mathcal{U}x)\d\fp(x),
\end{eqnarray*}
and this is exactly the statement of the Proposition. $\square$\\

For $L,T\in\mathbb{N}_{0}$ let us introduce
\begin{eqnarray*}
B_{L,T}(\pp)=\sum\limits_{n=0}^{\infty}\frac{1}{2^{n+1}\pp^{Tn}}\Big{(}\frac{\pp^{n}-1}{\pp^{n+1}-\pp^{n}}\Big{)}^{L}.
\end{eqnarray*}
For example,
\begin{eqnarray*}
B_{0,T}&=&\frac{\pp^{T}}{2\pp^{T}-1},\quad
B_{1,T}(\pp)=\frac{\pp^{T}}{(2\pp^{T}-1)(2\pp^{T+1}-1)},\\
B_{2,T}(\pp)&=&\frac{\pp^{T}(2\pp^{T+1}+1)}{(2\pp^{T+2}-1)(2\pp^{T+1}-1)(2\pp^{T}-1)},\\
B_{3,T}(\pp)&=&\frac{\pp^{T}(4\pp^{2T+3}+4\pp^{T+2}+4\pp^{T+1}+1)}{(2\pp^{T+3}-1)(2\pp^{T+2}-1)(2\pp^{T+1}-1)(2\pp^{T}-1)},\\
B_{4,T}(\pp)&=&\frac{\pp^{T}(2\pp^{T+2}+1)(4\pp^{2T+4}+6\pp^{T+3}+8\pp^{T+2}+6\pp^{T+1}+1)}{(2\pp^{T+4}-1)(2\pp^{T+3}-1)(2\pp^{T+2}-1)(2\pp^{T+1}-1)(2\pp^{T}-1)}.
\end{eqnarray*}
As it is easy to see, $B_{L,T}(\pp)$ are rational functions in $\pp$ for
$L,T\in\mathbb{N}_{0}$. Indeed,
\begin{eqnarray*}
B_{L,T}(\pp)=\frac{1}{(\pp-1)^{L}}\cdot\sum\limits_{n=0}^{\infty}\frac{1}{\pp^{Tn}2^{n+1}}
\Big{(}1-\frac{1}{\pp^{n}}\Big{)}^{L}=\frac{1}{2(\pp-1)^{L}}\cdot
\sum\limits_{s=0}^{L}(-1)^{s}\binom{L}{s}\sum\limits_{n=0}^{\infty}\frac{1}{2^{n}\pp^{n(s+T)}}=\\
\frac{\pp^{T}}{(\pp-1)^{L}}\cdot
\sum\limits_{s=0}^{L}(-1)^{s}\binom{L}{s}\frac{\pp^{s}}{2\pp^{s+T}-1}=
\frac{\pp^{T}\mathscr{R}_{L,T}(\pp)}{(2\pp^{T+L}-1)(2\pp^{T+L-1}-1)\cdot...\cdot(2\pp^{T+1}-1)(2\pp^{T}-1)},
\end{eqnarray*}
where $\mathscr{R}_{L,T}(\pp)$ are polynomials. This follows from the
observation
that $\pp=1$ is a root of numerator of multiplicity not less than $L$. \\

As in case $\pp=1$, our main concern are the moments of distributions $\fp(x)$,
which are defined by
\begin{eqnarray*}
m_{L}(\pp)&=&2\int\limits_{\I_{0}}x^{L}\d\fp(x)=
\int\limits_{\I_{\pp}}\Big{(}\frac{\pp x}{x+1}\Big{)}^{L}\d\fp(x)\\
&=&2\int\limits_{0}^{1}\mathfrak{X}^{L}(\pp,x)\d F(x)=
\lim_{n\rightarrow\infty}2^{2-n}\sum\limits_{a_{1}+a_{2}+...+a_{s}=n}[0,a_{1},a_{2},..,a_{s}]_{\pp}^{L}.,\\
M_{L}(\pp)&=&\int\limits_{\I_{\pp}}x^{L}\d\fp(x).
\end{eqnarray*}
Thus, if $\sup_{z\in\I_{\pp}}|z|=\rho_{\pp}>1$, which is finite for $\Re\pp\geq
1$, $\pp\neq 1$ (see Section 3), then $M_{L}(\pp)\leq \rho_{\pp}^{L}$.
\begin{prop}
The function $m_{L}(\pp)$ is analytic in the disc $|\pp-2|\leq 1$, including
its boundary. In particular, if in this disc
\begin{eqnarray*}
\widehat{m}_{L}(\pp):=\frac{m_{L}(\pp)}{\pp^{L}}=\sum\limits_{v=0}^{\infty}\eta_{v,L}(\pp-2)^{v},
\end{eqnarray*}
then for any $M\in\mathbb{N}$, one has the estimate $\eta_{v,L}\ll v^{-M}$ as
$v\rightarrow\infty$.
\label{prop4.2}
\end{prop}
{\it Proof. }The function $\mathfrak{X}(\pp,x)$ possesses a derivative in $\pp$
for $\Re\pp\geq 1$, $|\pp-2|\leq 1$, and these are bounded and continuous functions
for $x\in\mathbb{R}_{+}$. Therefore $m_{L}(\pp)$ has a derivative. For $\pp=1$,
there exists $\frac{\d^{M}}{\d\pp^{M}}\mathfrak{X}(\pp,x)\ll x^{M+1}$, and it
is a continuous function for irrational $x$. Additionally, $F'(x)=0$ for
$x\in\mathbb{Q}_{+}$. This proves the analyticity of $m_{L}(\pp)$ in the disc
$|\pp-2|\leq 1$. Then an estimate for the Taylor coefficients is the standard
fact from Fourier analysis. In fact,
\begin{eqnarray*}
\eta_{v,L}=\int\limits_{0}^{1}\widehat{m}_{L}(2+e^{2\pi i\vartheta})e^{-2\pi i
v\vartheta}\d \vartheta.
\end{eqnarray*}
The function $\widehat{m}_{L}(2+e^{2\pi i\vartheta})\in{\sf
C}^{\infty}(\mathbb{R})$, hence the iteration of integration by parts implies
the needed estimate. $\square$\\
\begin{prop}
Functions $M_{L}(\pp)$ and $m_{L}(\pp)$ are related via rational functions
$B_{L,T}(\pp)$ in the following way:
\begin{eqnarray*}
M_{L}(\pp)=\sum\limits_{s=0}^{L}m_{s}(\pp)B_{L-s,s}(\pp)\binom{L}{s}.
\end{eqnarray*}
\label{prop4.3}
\end{prop}
{\it Proof. }Indeed, this follows from the definitions and Proposition
\ref{prop4.1} in case $\omega(x)=x^{L}$. $\square$\\

Let us introduce, following \cite{ga1} in case $\pp=1$, the following
generating functions:
\begin{eqnarray}
\mathfrak{m}_{\pp}(t)&=&\sum\limits_{L=0}^{\infty}m_{L}(\pp)\frac{t^{L}}{L!}=
2\int\limits_{\I_{0}}e^{xt}\d\fp(x)=\int\limits_{\I_{\pp}}\exp\Big{(}\frac{\pp
xt}{x+1}\Big{)}\d\fp(x);\nonumber\\
G_{\pp}(z)&=&\sum\limits_{L=1}^{\infty}\frac{m_{L}(\pp)}{\pp^{L}}z^{L-1}=\int\limits_{\I_{\pp}}\frac{1}{x+1-z}\d
F_{\pp}(x)=\int\limits_{0}^{\infty}\frac{1}{\mathfrak{X}(\pp,x)+1-z}\d
F(x)\label{def}.
\end{eqnarray}
\indent The situation $\pp=2$ is particularly important, since all these
functions can be explicitly calculated, and it provides the case where all the
subsequent results can be checked directly and the starting point in proving
Theorem \ref{thm1.1}. Thus,
\begin{eqnarray*}
\mathfrak{m}_{2}(t)=e^{t},\quad G_{2}(z)=\frac{1}{2-z}.
\end{eqnarray*}
\indent By the definition, expressions $m_{L}(\pp)/\pp^{L}$ are Taylor
coefficients of $G_{\pp}(z)$ at $z=0$. Differentiation of $L-1$ times under the
integral defining $G_{\pp}(z)$, and substitution $z=1$ gives
\begin{eqnarray}
G_{\pp}^{(L-1)}(1)=(L-1)!\int\limits_{\I_{\pp}}\frac{1}{x^{L}}\d\fp(x)=(L-1)!M_{L}(\pp)\Rightarrow
G_{\pp}(z+1)=\sum\limits_{L=0}^{\infty}M_{L}(\pp)z^{L-1}\label{vienas},
\end{eqnarray}
with a radius of convergence equal to $\rho_{\pp}^{-1}$. As was proved in
\cite{ga1} and mentioned before, in case $\pp=1$ ($\rho_{1}=\infty$) this must
be interpreted that there exist all derivatives at $z=1$. The next Proposition
shows how symmetry property reflects in $\mathfrak{m}_{\pp}(t)$.
\begin{prop} One has
\begin{eqnarray*}
\mathfrak{m}_{\pp}(t)=e^{\pp t}\mathfrak{m}_{\pp}(-t).
\end{eqnarray*}
\label{prop4.4}
\end{prop}
{\it Proof. } Indeed,
\begin{eqnarray*}
\mathfrak{m}_{\pp}(t)=\int\limits_{\I_{\pp}}\exp\Big{(}\frac{\pp
xt}{x+1}\Big{)}\d\fp(x)=
\int\limits_{\I_{\pp}}\exp\Big{(}\pp t-\frac{\pp t}{x+1}\Big{)}\d\fp(x)=\\
e^{\pp t}\int\limits_{\I_{\pp}}\exp\Big{(}-\frac{\pp
t}{x+1}\Big{)}\d\fp(x)\mathop{=}^{x\rightarrow \frac{1}{x}}e^{\pp
t}\mathfrak{m}_{\pp}(-t).\quad\square
\end{eqnarray*}
This result allows to obtain linear relations among moments $m_{L}(\pp)$ and
the exact value of the first (trivial) moment $m_{1}(\pp)$.
\begin{cor} One has
\begin{eqnarray*}
m_{1}(\pp)=\frac{\pp}{2},\quad M_{1}(\pp)=\frac{\pp^{2}+2}{4\pp-2}.
\end{eqnarray*}
\label{cor4.5}
\end{cor}
{\it Proof. }Indeed, the last propositions implies
\begin{eqnarray*}
m_{L}(\pp)=\sum\limits_{s=0}^{L}\binom{L}{s}(-1)^{s}m_{s}(\pp)\pp^{L-s},\quad
L\geq 0.
\end{eqnarray*}
For $L=1$ this gives the first statement of the Corollary. Additionally,
Proposition \ref{prop4.3} for $L=1$ reads as
\begin{eqnarray*}
M_{1}(\pp)=\frac{\pp}{2\pp-1}\cdot m_{1}(\pp)+\frac{1}{2\pp-1},
\end{eqnarray*}
and we are done. $\square$
\section{Three term functional equation}
\begin{thm} The function $G_{\pp}(z)$ can be extended to analytic function in the
domain $\mathbb{C}\setminus (\I_{\pp}+1)$. It satisfies the functional equation
\begin{eqnarray}
\frac{1}{z}+\frac{\pp}{z^{2}}G_{\pp}\Big{(}\frac{\pp}{z}\Big{)}+2G_{\pp}(z+1)=\pp
G_{\pp}(\pp z), \text{ for }z\notin\frac{\I_{\pp}+1}{\pp}. \label{fep}
\end{eqnarray}
Its consequence is the symmetry property
\begin{eqnarray*}
G_{\pp}(z+1)=-\frac{1}{z^{2}}G_{\pp}\Big{(}\frac{1}{z}+1\Big{)}-\frac{1}{z}.
\end{eqnarray*}
Moreover, $G_{\pp}(z)\rightarrow 0$ if
$\textrm{dist}(z,\I_{\pp})\rightarrow\infty$.\\
Conversely - the function satisfying this functional equation and regularity
property is unique.
\label{thm5.1}
\end{thm}
{\it Proof. }Let $w(x,z)=\frac{1}{x+1-z}$. Then it is straightforward to check
that
\begin{eqnarray*}
w(\frac{x+1}{\pp},z+1)=\pp\cdot w(x,\pp z),\\
w(\frac{\pp}{x+1},z+1)=-\frac{\pp}{z^{2}}w(x,\frac{\pp}{z})-\frac{1}{z}.
\end{eqnarray*}
Thus, for $|\pp-2|\leq 1$, $\pp\neq 2$,
\begin{eqnarray*}
2G_{\pp}(z+1)&=&2\int\limits_{\I_{0}}w(x,z+1)\d\fp(x)+
2\int\limits_{\I_{\pp}\setminus\I_{0}}w(x,z+1)\d\fp(x)\\
&=&2\int\limits_{\I_{\pp}}w(\frac{\pp x}{x+1},z+1)\d\fp\Big{(}\frac{\pp
x}{x+1}\Big{)}+2\int\limits_{\I_{\pp}}w(\frac{x+1}{\pp},z+1)\d\fp\Big{(}\frac{x+1}{\pp}\Big{)}\\
&=&\int\limits_{\I_{\pp}}w(\frac{\pp}{x+1},z+1)\d\fp(x)+\int\limits_{\I_{\pp}}w(\frac{x+1}{\pp},z+1)\d\fp(x)\\
&=&-\frac{1}{z}-\frac{\pp}{z^{2}}G_{\pp}\Big{(}\frac{\pp}{z}\Big{)}+\pp
G_{\pp}(\pp z).
\end{eqnarray*}
(In the first integral we used a substitution $x\rightarrow\frac{1}{x}$). The
functional equation holds in case $\pp=2$ as well, which can be checked
directly. The holomorphicity of $G_{\pp}(z)$ follows exactly as in case
$\pp=1$ \cite{ga1}. All we need is the first integral in (\ref{def}) and
the fact that $\I_{\pp}$ is a closed set.\\
\indent As was mentioned, the uniqueness of a function satisfying (\ref{fep})
for $\pp=1$ was proved in \cite{ga1}. Thus, the converse implication follows
from analytic continuation principle for the function in two complex variables
$(\pp,z)$ (see Lemma \ref{lem6.2} below, where the proof in case $\pp=2$ is
presented. Similar argument works for general $\pp$). $\square$

\begin{cor}Let $\pp\neq 1$, and $\mathscr{C}$ be any closed smooth contour which rounds the
curve $\I_{\pp}+1$ once in the positive direction. Then
\begin{eqnarray*}
\frac{1}{2\pi i}\oint\limits_{\mathscr{C}}G_{\pp}(z)\d z=-1.
\end{eqnarray*}
\label{cor5.2}
\end{cor}
{\it Proof. }Indeed, this follows from the functional equation (\ref{fep}), as well as from
the symmetry property. It is enough to take a sufficiently large circle
$\mathscr{C}=\{|z|=R\}$ such that $\mathscr{C}^{-1}+1$ is contained in a small
neighborhood of $z=1$, for which
$(\mathscr{C}^{-1}+1)\cap(\I_{\pp}+1)=\emptyset$. This is possible since
$0\notin\I_{\pp}$ (see Theorem \ref{thm3.1}). $\square$\\

We finish with providing an integral equation for $\mathfrak{m}_{\pp}(t)$. We
indulge in being concise since the argument directly generalizes the one used
in \cite{ga1} to prove the integral functional equation for $\mathfrak{m}(t)$
(in our notation, this is $\mathfrak{m}_{1}(t)$).
\begin{prop}
Let $1\leq\pp<\infty$ be real. Then the function $\mathfrak{m}_{\pp}(t)$
satisfies the boundary condition $\mathfrak{m}_{\pp}(0)=1$, regularity property
$\mathfrak{m}_{\pp}(-t)\ll e^{-\sqrt{t\log 2}}$, and the integral equation
\begin{eqnarray*}
\mathfrak{m}_{\pp}(-s)=\int\limits_{0}^{\infty}\mathfrak{m}'_{\pp}(-t)\Big{(}2e^{s}
J_{0}(2\sqrt{\pp st})-J_{0}(2\sqrt{st})\Big{)}\d t,\quad s\in\mathbb{R}_{+}.
\end{eqnarray*}
\label{prop5.3}
\end{prop}
For instance, in the case $\pp=1$ this reduces to (\ref{char}), and in the case $\pp=2$ this reads as
\begin{eqnarray*}
2e^{s}\int\limits_{0}^{\infty}e^{-t}J_{0}(2\sqrt{2st})\d
t=2e^{s}e^{-2s}=e^{-s}+e^{-s}=e^{-s}+\int\limits_{0}^{\infty}e^{-t}J_{0}(2\sqrt{st})\d
t,
\end{eqnarray*}
which is an identity \cite{wats}.\\
{\it Proof. }Indeed, the functional equation for $G_{\pp}(z)$ in the region $\Re z<-1$
in terms of $\mathfrak{m}'_{\pp}(t)$ reads as
\begin{eqnarray*}
\frac{1}{z}=\int\limits_{0}^{\infty}\mathfrak{m}'_{\pp}(-t)
\Big{(}\frac{2}{z+1}e^{\frac{\pp t}{z+1}}+\frac{1}{z}e^{tz}-
\frac{1}{z}e^{\frac{t}{z}}\Big{)}\d t.
\end{eqnarray*}
Now, multiply this by $e^{-sz}$ and integrate over $\Re z=-\sigma<-1$,
where $s>0$ is real. All the remaining steps are exactly the same as in \cite{ga1}. $\square$\\

\noindent {\it Remark. }If $\pp\neq 1$, the regularity bound is easier than in case
$\pp=1$. Take, for example, $1<\pp< 2$. Then
\begin{eqnarray*}
|\mathfrak{m}_{\pp}(t)|\leq\int\limits_{\pp-1}^{\frac{1}{\pp-1}}\Big{|}\exp\Big{(}\frac{\pp
xt}{x+1}\Big{)}\Big{|}\d\fp
(x)<\int\limits_{\pp-1}^{\frac{1}{\pp-1}}e^{t}\d\fp(x)=e^{t}.
\end{eqnarray*}
Thus, Proposition \ref{prop4.4} gives $|\mathfrak{m}_{\pp}(-t)|<e^{(1-\pp)t}$.
The same argument shows that for $\pp>2$ we have
$|\mathfrak{m}_{\pp}(-t)|<e^{-t}$.

\section{The proof: approach through $\pp=2$}
Let us rewrite the functional equation for $G_{\pp}(z)=G(\pp,z)$ as
\begin{eqnarray}
\frac{1}{z}+\frac{\pp}{z^{2}}G\Big{(}\pp,\frac{\pp}{z}\Big{)}+2G(\pp,z+1)=\pp
G(\pp,\pp z).\label{dvi}
\end{eqnarray}
Direct induction shows that the following ``chain-rule" holds:
\begin{eqnarray}
\frac{\p^{n}}{\p\pp^{n}}\Big{(}\pp G(\pp,\pp z)\Big{)}=
\sum\limits_{i+j=n}\binom{n}{j}\pp\frac{\p^{i}\p^{j}}{\p\pp^{i}\p
z^{j}}G(\pp,\pp z)z^{j}+\nonumber\\
\sum\limits_{i+j=n-1}n\binom{n-1}{j}\frac{\p^{i}\p^{j}}{\p\pp^{i}\p
z^{j}}G(\pp,\pp z)z^{j},\label{dal}
\end{eqnarray}
where in the summation it is assumed that $i,j\geq 0$.\\
\indent Now we will provide rigorous calculations which yield explicit series for
$G(\pp,z)$ in terms of powers of $(\pp-2)$ and certain rational functions. The
function $G(\pp,z)$ is analytic in $\{|\pp-2|\leq 1\}\times \{|z|\leq \frac{3}{4}\}$.
This follows from Theorem \ref{thm3.1} and integral representation (\ref{def}).
Thus, for $\{|\pp-2|< 1\}\times \{|z|\leq\frac{3}{4}\}$ it has a
Taylor expansion
\begin{eqnarray}
G(\pp,z)=\sum\limits_{L=1}^{\infty}\sum\limits_{v=0}^{\infty}\eta_{v,L}\cdot
z^{L-1}(\pp-2)^{v}.\label{taylor}
\end{eqnarray}
Moreover, the function $G(2+e^{2\pi i\vartheta},\frac{3}{4}e^{2\pi i\varphi})\in{\sf
C}^{\infty}(\mathbb{R}\times\mathbb{R})$, and it is double-periodic. Thus,
\begin{eqnarray*}
\eta_{v,L}=\Big{(}\frac{4}{3}\Big{)}^{L-1}\int\limits_{0}^{1}\int\limits_{0}^{1}G(2+e^{2\pi
i\vartheta},\frac{3}{4}e^{2\pi i\varphi})e^{-2\pi i v\vartheta-2\pi i
(L-1)\varphi}\d\vartheta\d\varphi,\quad v\geq 0,\quad L\geq 1.
\end{eqnarray*}
A standard trick from Fourier analysis (using iteration of integration by
parts) shows that $\eta_{v,L}\ll_{M} (4/3)^{L}\cdot(Lv)^{-M}$ for any $M\in\mathbb{N}$. Thus,
(\ref{taylor}) holds for $(\pp,z)\in\{|\pp-2|\leq 1\}\times \{|z|\leq 3/4\}$.\\
\indent Our idea is a simple one. Indeed, let us look at (\ref{def}). This
implies the Taylor series for
$m_{L}(\pp)/\pp^{L}=\sum_{v=0}^{\infty}\eta_{v,L}(\pp-2)^{v}$, convergent in
the disc $|\pp-2|\leq 1$. Due to the absolute convergence, the order of
summation in (\ref{taylor}) is not essential. This yields
\begin{eqnarray*}
G(\pp,z)=\sum\limits_{v=0}^{\infty}(\pp-2)^{v}\Big{(}\sum\limits_{L=1}^{\infty}\eta_{v,L}\cdot
z^{L-1}\Big{)}.
\end{eqnarray*}
Therefore, let
\begin{eqnarray*}
\frac{1}{n!}\frac{\p^{n}}{\p\pp^{n}}G(\pp,z)\Big{|}_{\pp=2}=\mathbf{H}_{n}(z)=\sum\limits_{L=1}^{\infty}\eta_{n,L}\cdot
z^{L-1}.
\end{eqnarray*}
We already know that $\mathbf{H}_{0}(z)=\frac{1}{2-z}$. Though $m_{L}(\pp)$ are
obviously highly transcendental (and mysterious) functions, the series for $\mathbf{H}_{n}(z)$
is in fact a rational function in $z$, and this is the main point of our
approach. Moreover, we will show that
\begin{eqnarray*}
\mathbf{H}_{n}(z)=\frac{\mathscr{B}_{n}(z)}{(z-2)^{n+1}},
\end{eqnarray*}
where $\mathscr{B}_{n}(z)$ is a polynomial with rational coefficients of degree
$n-1$ with the reciprocity property
$\mathscr{B}_{n}(z+1)=(-1)^{n}z^{n-1}\mathscr{B}_{n}(\frac{1}{z}+1)$,
$\mathscr{B}_{n}(0)=0$. We argue by induction on $n$.
First we need an auxiliary lemma.\\

Let $\mathbb{Q}[z]_{n-1}$ denote the linear space of dimension $n$ of
polynomials of degree $\leq n-1$ with rational coefficients. Consider a
following linear map
$\mathcal{L}_{n-1}:\mathbb{Q}[z]_{n-1}\rightarrow\mathbb{Q}[z]_{n-1}$, defined
by
\begin{eqnarray*}
\mathcal{L}_{n-1}(P)(z)=P(z+1)-\frac{1}{2^{n+1}}P(2z)+\frac{(-1)^{n+1}}{2^{n+1}}P\Big{(}\frac{2}{z}\Big{)}z^{n-1}.
\end{eqnarray*}
\begin{lem}
$\text{det}(\mathcal{L}_{n-1})\neq 0$. Accordingly, $\mathcal{L}_{n-1}$ is the
isomorphism.
\label{lem6.1}
\end{lem}
\noindent {\it Remark. }Let $m=\big{[}\frac{n}{2}\big{]}$. Then it can be proved that
indeed $\text{det}(\mathcal{L}_{n-1})=\frac{\prod_{i=1}^{m}(4^i-1)}{2^{m^2+m}}$.\\

{\it Proof. }Suppose $P\in\text{ker}(\mathcal{L}_{n-1})$. Then a rational
function $\mathbf{H}(z)=\frac{P(z)}{(z-2)^{n+1}}$ satisfies the three term
functional equation
\begin{eqnarray}
\mathbf{H}(z+1)-\mathbf{H}(2z)+\mathbf{H}\Big{(}\frac{2}{z}\Big{)}\frac{1}{z^{2}}=0,\quad
z\neq 1.\label{hom}
\end{eqnarray}
Also, $\mathbf{H}(z)=o(1)$, as $z\rightarrow\infty$. Now the result follows
from the next
\begin{lem}
Let $\Upsilon(z)$ be any analytic function in the domain
$\mathbb{C}\setminus\{1\}$. Then if $\mathbf{H}(z)$ is a solution of the
equation
\begin{eqnarray*}
\mathbf{H}(z+1)-\mathbf{H}(2z)+\mathbf{H}\Big{(}\frac{2}{z}\Big{)}\frac{1}{z^{2}}=\Upsilon(z),
\end{eqnarray*}
such that $\mathbf{H}(z)\rightarrow 0$ as $z\rightarrow\infty$, $\mathbf{H}(z)$
is analytic in $\mathbb{C}\setminus\{2\}$, then such $\mathbf{H}(z)$ is unique.
\label{lem6.2}
\end{lem}
{\it Proof. }All we need is to show that with the imposed diminishing
condition, homogeneous equation (\ref{hom}) admits only the solution
$\mathbf{H}(z)\equiv 0$. Indeed, let $\mathbf{H}(z)$ be such a solution. Put
$z\rightarrow 2^{n}z+1$. Thus,
\begin{eqnarray*}
\mathbf{H}(2^{n}z+2)-\mathbf{H}(2^{n+1}z+2)+\frac{1}{(2^{n}z+1)^{2}}
\mathbf{H}\Big{(}\frac{2}{2^{n}z+1}\Big{)}=0.
\end{eqnarray*}
 This is valid for $z\neq 0$ (since $\mathbf{H}(z)$ is allowed
to have a singularity at $z=2$). Now sum this over $n\geq 0$. Due to the
diminishing assumption, one gets (after additional substitution $z\rightarrow
z-2$)
\begin{eqnarray*}
\mathbf{H}(z)=-\sum\limits_{n=0}^{\infty}\frac{1}{(2^{n}z-2^{n+1}+1)^2}\mathbf{H}\Big{(}\frac{2}{2^{n}z-2^{n+1}+1}\Big{)}.
\end{eqnarray*}
For clarity, put $z\rightarrow -z$ and consider a function
$\widehat{\mathbf{H}}(z)=\mathbf{H}(-z)$. Thus,
\begin{eqnarray*}
\widehat{\mathbf{H}}(z)=-\sum\limits_{n=0}^{\infty}\frac{1}{(2^{n}z+2^{n+1}-1)^2}\widehat{\mathbf{H}}\Big{(}\frac{2}{2^{n}z+2^{n+1}-1}\Big{)}.
\end{eqnarray*}
Consider this for $z\in[0,2]$. As can be easily seen, then all arguments on the
right also belong to this interval. We want to prove the needed result simply
by applying the maximum argument. The last identity is still insufficient. For
this reason consider its second iteration. This produces a series
\begin{eqnarray*}
\widehat{\mathbf{H}}(z)=\sum\limits_{n,m=0}^{\infty}
\frac{1}{(2^{n+m+1}z+2^{n+m+2}-2^{n}z-2^{n+1}+1)^2}\widehat{\mathbf{H}}\Big{(}\omega_{m}\circ\omega_{n}(z)\Big{)},
\end{eqnarray*}
where $\omega_{n}(z)=\frac{2}{2^{n}z+2^{n+1}-1}$. As said,
$\omega_{m}\circ\omega_{n}(z)\in[0,2]$ for $z\in[0,2]$. Since a function
$\widehat{\mathbf{H}}(z)$ is continuous in the interval $[0,2]$, let $z_{0}\in[0,2]$ be
such that
$M=|\widehat{\mathbf{H}}(z_{0})|=\sup_{z\in[0,2]}|\widehat{\mathbf{H}}(z)|$.
Consider the above expression for $z=z_{0}$. Thus,
\begin{eqnarray*}
M=|\widehat{\mathbf{H}}(z_{0})|\leq
\sum\limits_{n,m=0}^{\infty}\Big{|}\frac{1}{(2^{n+m+1}z_{0}+2^{n+m+2}-2^{n}z_{0}-2^{n+1}+1)^2}
\widehat{\mathbf{H}}\Big{(}\omega_{m}\circ\omega_{n}(z_{0})\Big{)}\Big{|}\leq\\
M\sum\limits_{n,m=0}^{\infty}\frac{1}{(2^{n+m+2}-2^{n+1}+1)^{2}}=0.20453_{+}M.
\end{eqnarray*}
This is contradictory unless $M=0$. By the principle of analytic continuation,
$\mathbf{H}(z)\equiv 0$, and this proves the Lemma. $\square$\\

\noindent{\it Remark. } Direct inspection of the proof reveals that the statement of
Lemma still holds with a weaker assumption that $\mathbf{H}(z)$ is
real-analytic function on $(-\infty,0]$.\\

Now, let us differentiate (\ref{dvi}) $n$ times with respect to $\pp$, use
(\ref{dal}) and afterwards substitute $\pp=2$. This gives
\begin{eqnarray}
\sum\limits_{j=1}^{n}\frac{2}{j!}\frac{\p^{j}}{\p
z^{j}}\mathbf{H}_{n-j}(2z)z^{j}+
\sum\limits_{j=0}^{n-1}\frac{1}{j!}\frac{\p^{j}}{\p
z^{j}}\mathbf{H}_{n-j-1}(2z)z^{j}-\nonumber\\
\sum\limits_{j=1}^{n}\frac{2}{j!}\frac{\p^{j}}{\p
z^{j}}\mathbf{H}_{n-j}\Big{(}\frac{2}{z}\Big{)}\frac{1}{{z}^{j+2}}-
\sum\limits_{j=0}^{n-1}\frac{1}{j!}\frac{\p^{j}}{\p
z^{j}}\mathbf{H}_{n-j-1}\Big{(}\frac{2}{z}\Big{)}\frac{1}{z^{j+2}}=\nonumber\\
2\mathbf{H}_{n}(z+1)-2\mathbf{H}_{n}(2z)+2\mathbf{H}_{n}\Big{(}\frac{2}{z}\Big{)}\frac{1}{z^{2}}.
\label{lygtis}
\end{eqnarray}
We note that this implies the reciprocity property
\begin{eqnarray*}
\mathbf{H}_{n}(z+1)=-\frac{1}{z^{2}}\mathbf{H}_{n}\Big{(}\frac{1}{z}+1\Big{)},
\quad n\geq 1.
\end{eqnarray*}
{\it A posteriori}, this clarifies how the identity $F(x)+F(1/x)=1$ reflects in
the series for $G(z)$, as stated in Theorem \ref{thm1.1}: reciprocity property
(non-homogeneous for $n=0$ and homogeneous for $n\geq 1$) is reflected in each
of the summands separately, whereas the three term functional equation heavily
depends on inter-relations among $\mathbf{H}_{n}(z)$.\\

Now, suppose we know all $\mathbf{H}_{j}(z)$ for $j<n$.
\begin{lem}
The left hand side of the equation (\ref{lygtis}) is of the form
\begin{eqnarray*}
\text{l.h.s.}=\frac{\mathscr{J}_{n}(z)}{(z-1)^{n+1}},
\end{eqnarray*}
where $\mathscr{J}_{n}(z)\in\mathbb{Q}[z]_{n-1}$.
\label{lem6.3}
\end{lem}
{\it Proof. } First, as it is clear from the appearance of l.h.s., we need to
verify that $z$ does not divide a denominator, if l.h.s. is represented as a
quotient of two co-prime polynomials. Indeed, using the symmetry property in
(\ref{dvi}) for the term $G(\pp,\frac{\pp}{z})$, we obtain the three term
functional equation of the form
\begin{eqnarray*}
-\frac{1}{\pp-z}-\frac{\pp}{(\pp-z)^{2}}G\Big{(}\pp,\frac{\pp}{\pp-z}\Big{)}+2G(\pp,z+1)=\pp
G(\pp,\pp z).
\end{eqnarray*}
Let us perform the same procedure which we followed to arrive at the equation
(\ref{lygtis}). Thus, differentiation $n$ times with respect to $\pp$ and
substitution $\pp=2$ gives the expression of the form
\begin{eqnarray*}
\text{l.h.s.}_{2}=2\mathbf{H}_{n}(z+1)-2\mathbf{H}_{n}(2z)-2\mathbf{H}_{n}\Big{(}\frac{2}{2-z}\Big{)}\frac{1}{(2-z)^{2}},
\end{eqnarray*}
where $\text{lh.s.}_{2}$ is expressed in terms of $\mathbf{H}_{j}(z)$ for
$j<n$. Nevertheless, this time the common denominator of $\text{l.h.s.}_{2}$ is
of the form $(z-1)^{n+1}(z-2)^{n+2}$. As a corollary, $z$ does not divide it.
Finally, due to the reciprocity property, for $n\geq 1$ one has
\begin{eqnarray*}
\mathbf{H}_{n}\Big{(}\frac{2}{2-z}\Big{)}\frac{1}{(2-z)^{2}}=-\mathbf{H}_{n}\Big{(}\frac{2}{z}\Big{)}\frac{1}{z^{2}}.
\end{eqnarray*}
This shows that actually $\text{l.h.s.}=\text{l.h.s.}_{2}$, and therefore if
this is expressed as a quotient of two polynomials in lowest terms, the
denominator is a power of $(z-1)$. Finally, it is obvious that this exponent is
exactly $n+1$, and one easily verifies that $\deg \mathscr{J}_{n}(z)\leq n-1$.
(Possibly, $\mathscr{J}_{n}(z)$ can be divisible by $(z-1)$, but this does not
have an
impact on the result). $\square$\\

{\bf Proof of Theorem \ref{thm1.1}}. Now, using Lemma \ref{lem6.1}, we
inherit that there exists a unique polynomial $\mathscr{B}_{n}(z)$ of degree
$\leq n-1$ such that
$\mathscr{B}_{n}(z)=\frac{1}{2}\mathcal{L}_{n-1}^{-1}(\mathscr{J}_{n})(z)$.
Summarizing, $\mathbf{H}_{n}(z)=\frac{\mathscr{B}_{n}(z)}{(z-2)^{n+1}}$ solves
the equation (\ref{lygtis}). On the other hand, Lemma \ref{lem6.2} implies that
the solution of (\ref{lygtis}) we obtained is indeed the unique one. This
reasoning proves that for $|\pp-2|\leq 1$, $|z|\leq 1$ we have the series
\begin{eqnarray*}
G(\pp,z)=\sum\limits_{n=0}^{\infty}(\pp-2)^{n}\cdot\mathbf{H}_{n}(z).
\end{eqnarray*}
This finally establishes the validity of Theorem \ref{thm1.1}. Note also that each summand
satisfies the symmetry property. The series converges absolutely for any $z$,
$|z|\leq3/4$, and if this holds for $z$, the same does hold for $\frac{z}{z-1}$,
which gives the circle $|z+9/7|\leq 12/7$. $\square$\\

Curiously, one could formally verify that the function defined by this series
does indeed satisfy (\ref{fep}). Indeed, using (\ref{lygtis}), we get:
\begin{eqnarray*}
2G(\pp,z+1)=2\mathbf{H}_{0}(z+1)+2\sum\limits_{n=1}^{\infty}(\pp-2)^{n}\mathbf{H}_{n}(z+1)=\\
2\mathbf{H}_{0}(z+1)+\sum\limits_{n=1}^{\infty}(\pp-2)^{n}\Bigg{(}\sum\limits_{j=0}^{n}\frac{2}{j!}\frac{\p^{j}}{\p
z^{j}}\mathbf{H}_{n-j}(2z)z^{j}+
\sum\limits_{j=0}^{n-1}\frac{1}{j!}\frac{\p^{j}}{\p
z^{j}}\mathbf{H}_{n-j-1}(2z)z^{j}-\\
\sum\limits_{j=0}^{n}\frac{2}{j!}\frac{\p^{j}}{\p
z^{j}}\mathbf{H}_{n-j}\Big{(}\frac{2}{z}\Big{)}\frac{1}{{z}^{j+2}}-
\sum\limits_{j=0}^{n-1}\frac{1}{j!}\frac{\p^{j}}{\p
z^{j}}\mathbf{H}_{n-j-1}\Big{(}\frac{2}{z}\Big{)}\frac{1}{z^{j+2}}\Bigg{)}
\end{eqnarray*}
Denote $n-j=s$. Then interchanging the order of summation for the first term of
the sum in the brackets, we obtain:
\begin{eqnarray*}
2\sum\limits_{n=1}^{\infty}(\pp-2)^{n}\sum\limits_{j=0}^{n}\frac{1}{j!}\frac{\p^{j}}{\p
z^{j}}\mathbf{H}_{n-j}(2z)z^{j}=2\sum\limits_{s=0}^{\infty}\sum\limits_{j=0}^{\infty}
(\pp-2)^{j+s}\frac{1}{j!}\frac{\p^{j}}{\p
z^{j}}\mathbf{H}_{s}(2z)z^{j}-2\mathbf{H}_{0}(2z)=\\
2\sum\limits_{s=0}^{\infty}(\pp-2)^{s}\mathbf{H}_{s}(2z+(\pp-2)z)-2\mathbf{H}_{0}(2z)=2G(\pp,\pp z)-2\mathbf{H}_{0}(2z).
\end{eqnarray*}
The same works for the second sum:
\begin{eqnarray*}
\sum\limits_{n=1}^{\infty}(\pp-2)^{n}\sum\limits_{j=0}^{n-1}\frac{1}{j!}\frac{\p^{j}}{\p
z^{j}}\mathbf{H}_{n-j-1}(2z)z^{j}=(\pp-2)G(\pp,\pp z).
\end{eqnarray*}
We perform the same interchange of summation for the second and the third
summand respectively. Thus, this yields
\begin{eqnarray*}
2G(\pp,z+1)=\pp G(\pp,\pp z)-\frac{\pp}{z^{2}}G\Big{(}\pp, \frac{\pp}{z}\Big{)}+2\mathbf{H}_{0}(z+1)-
2\mathbf{H}_{0}(2z)+\frac{2}{z^{2}}\mathbf{H}_{0}\Big{(}\frac{2}{z}\Big{)}=\\
\pp G(\pp,\pp z)-\frac{\pp}{z^{2}}G\Big{(}\pp, \frac{\pp}{z}\Big{)}-\frac{1}{z}.
\end{eqnarray*}
On the other hand, it is unclear how one can make this argument to work. This
would require rather detailed investigation of the linear map
$\mathcal{L}_{n-1}$ and recurrence (\ref{lygtis}), and this seems to be
very technical.

\appendix
\section{}
\subsection{Approach through $\pp=0$}
With a slight abuse of notation, we will use the expression
$\frac{\p^{s}}{\p\pp^{s}}G(0,z)$ to denote
$\frac{\p^{s}}{\p\pp^{s}}G(\pp,z)\big{|}_{\pp=0}$ for $s\in\mathbb{N}_{0}$.
Though the function $G(\pp,z)$ is defined only for $\Re \pp\geq 1$,
$z\notin(\I_{\pp}+1)$, assume that we are able to prove that it is analytic in
$\pp$ in a certain wider domain containing a disc $|\pp|<\varpi$, $\varpi>0$.
These are only formal calculations, but they unexpectedly yield series (\ref{exr}) (see
Section 1), and numerical calculations do strongly confirm the validity of
it.\\
\indent Thus, substitution $\pp=0$ into (\ref{dvi}) gives $G(0,z)=\frac{1}{2(1-z)}$.
Partial differentiation of (\ref{dvi}) with respect to $\pp$, and
consequent substitution $\pp=0$ gives
\begin{eqnarray*}
\frac{1}{z^{2}}G(0,0)+2\frac{\p}{\p\pp}G(0,z+1)=G(0,0)\Rightarrow
\frac{\p}{\p\pp}G(0,z)=\frac{(z-1)^{2}-1}{4(z-1)^{2}}.
\end{eqnarray*}
In the same fashion, differentiating the second time, we obtain
$\frac{\p^{2}}{\p\pp^{2}}G(0,z)=\frac{(z-1)^{4}-1}{2(z-1)^{3}}$.  In general, differentiating
(\ref{dvi}) $n\geq 1$ times with respect to $\pp$, using (\ref{dal}), and substituting $\pp=0$, we
obtain:
\begin{eqnarray*}
2\frac{\p^{n}}{\p\pp^{n}}G(0,z+1)=\sum\limits_{i+j=n-1}n\binom{n-1}{j}\frac{\p^{i}\p^{j}}{\p\pp^{i}\p
z^{j}}G(0,0)\Big{(}z^{j}-\frac{1}{z^{j+2}}\Big{)}.
\end{eqnarray*}
Let
\begin{eqnarray*}
\frac{1}{n!}\cdot\frac{\p^{n}}{\p\pp^{n}}G(0,z)=\overline{\mathbf{Q}}_{n}(z).
\end{eqnarray*}
Then
\begin{eqnarray*}
2\overline{\mathbf{Q}}_{n}(z+1)=\sum\limits_{j=0}^{n-1}\frac{1}{j!}\frac{\p^{j}}{\p
z^{j}}\overline{\mathbf{Q}}_{n-j-1}(0)\Big{(}z^{j}-\frac{1}{z^{j+2}}\Big{)}.
\end{eqnarray*}
Consequently, we have a recurrent formula to compute rational functions
$\overline{\mathbf{Q}}(z)$. Let
$\mathbf{Q}_{n}(z)=\overline{\mathbf{Q}}_{n}(z+1)$. Thus,
\begin{eqnarray*}
\mathbf{Q}_{n}(z)=\frac{(z+1)(z-1)\mathscr{D}_{n}(z)}{z^{n+1}},\quad n\geq 1,
\end{eqnarray*}
where $\mathscr{D}_{n}$ are polynomials of degree $2n-2$ with the reciprocity
property $\mathscr{D}_{n}(z)=z^{2n-2}\mathscr{D}_{n}\Big{(}\frac{1}{z}\Big{)}$
(this is obvious from the recurrence relation which defines
$\mathbf{Q}_{n}(z)$). Moreover, the coefficients of $\mathscr{D}_{n}$ are
$\mathbb{Q}_{p}$ integers for any prime $p\neq 2$. These calculations yield a
following formal result:
\begin{eqnarray*}
G(\pp,z)``="\sum\limits_{n=0}^{\infty}\pp^{n}\cdot\mathbf{Q}_{n}(z-1)=
\sum\limits_{n=0}^{\infty}\pp^{n}\frac{z(z-2)\mathscr{D}_{n}(z-1)}{(z-1)^{n+1}}.
\end{eqnarray*}
This produces the ``series" for the second and higher moments of the form
\begin{eqnarray*}
m_{2}(\pp)=\pp^{2}\cdot\sum\limits_{n=0}^{\infty}\pp^{n}\mathbf{Q}'_{n}(-1).
\end{eqnarray*}
\indent In particular, inspection of the table in Section 1 (where the initial
values for $\mathbf{Q}'_{n}(-1)$ are listed) shows that this series for $\pp=1$
does not converge. However, the Borel sum is properly defined and it converges
exactly to the value $m_{2}$. This gives empirical evidence for the validity of
(\ref{exr}). The principles of Borel summation also suggest the mysterious fact
that indeed $G(\pp,z)$ analytically extends to the
interval $\pp\in[0,1]$.\\

\indent Additionally, numerical calculations reveal the following fact: the
sequence $\sqrt[n]{|\mathbf{Q}'_{n}(-1)|}$ is monotonically increasing
(apparently, tends to $\infty$), while
$\frac{1}{n}\log|\mathbf{Q}'_{n}(-1)|-\log n$ monotonically decreases
(apparently, tends to $-\infty$). Thus,
\begin{eqnarray*}
A^{n}<|\mathbf{Q}'_{n}(-1)|<(cn)^{n},
\end{eqnarray*}
for $c=0.02372$ and $A=3.527$, $n\geq 150$.  We do not have enough evidence to
conjecture the real growth of this sequence. If $c=c(n)\rightarrow 0$, as
$n\rightarrow\infty$, then the function
\begin{eqnarray*}
\Lambda(t)=\sum\limits_{n=0}^{\infty}\frac{\mathbf{Q}'_{n}(-1)}{n!}t^{n}
\end{eqnarray*}
is entire, and in case $L=2$, result (\ref{exr}) is equivalent to the fact that
\begin{eqnarray*}
\int\limits_{0}^{\infty}\Lambda(t)e^{-t}\d t=m_{2}.
\end{eqnarray*}
\subsection{Auxiliary lemmas}
These lemmas are needed in Section 3. For $a,b\in\mathbb{N}$, $\pp\in\mathbb{C}$, $|\pp-2|\leq 1$, define rational functions
\begin{eqnarray*}
W_{a}(\pp)=\frac{\pp^{a}-1}{\pp^{a+1}-\pp^{a}},\quad T_{a,b}(\pp)=W^{-1}_{a}(\pp)W^{-1}_{b}(\pp)\pp^{-a}=\frac{(\pp-1)^{2}\pp^{b}}{(\pp^{a}-1)(\pp^{b}-1)}.
\end{eqnarray*}
Let us define constants
\begin{eqnarray*}
\mu(a,b)=\sup\limits_{\pp\in\mathbb{C},|\pp-2|\leq 1}|T_{a,b}(\pp)|-\Re(T_{a,b}(\pp)).
\end{eqnarray*}
The following table provides some initial values for constants $\mu(a,b)$, computed numerically.
\begin{center}
\begin{tabular}{|r| r| r| r| r| r| r|}
\hline
$b\setminus a$   & $1$ & $2$ & $3$ & $4$& $5$ & $6$\\
\hline
1 &$0.25000000$ &$0.01250000$ &$0.00780868$ &$0.03343231$ &$0.05778002$&$0.07712952$\\
2 &$0.29846114$ &$0.03125000$ &$0.00159908$ &$0.01212467$ &$0.02539758$&$0.03645721$\\
3 &$0.35999295$ &$0.05097235$ &$0.00647895$ &$0.00676996$ &$0.01624300$&$0.02437494$\\
4 &$0.41433340$ &$0.07007201$ &$0.01316542$ &$0.00500146$ &$0.01287728$&$0.01963810$\\
5 &$0.45590757$ &$0.08747624$ &$0.02069451$ &$0.00437252$ &$0.01163446$&$0.01781467$\\
6 &$0.48390408$ &$0.10255189$ &$0.02845424$ &$0.00812804$ &$0.01125132$&$0.01728395$\\
7 &$0.49985799$ &$0.11503743$ &$0.03601828$ &$0.01200557$ &$0.01120308$&$0.01729854$\\
8 &$0.50642035$ &$0.12494927$ &$0.04309384$ &$0.01611126$ &$0.01125789$&$0.01748823$\\
9 &$0.50681483$ &$0.13248892$ &$0.04949922$ &$0.02025219$ &$0.01132055$&$0.01767914$\\
10&$0.50452450$ &$0.13796512$ &$0.05514483$ &$0.02427779$ &$0.01136245$&$0.01780892$\\
11&$0.50218322$ &$0.14173414$ &$0.06001269$ &$0.02807992$ &$0.01138335$&$0.01787452$\\
12&$0.50070286$ &$0.14415527$ &$0.06413550$ &$0.03158969$ &$0.01139099$&$0.01789618$\\
13&$0.49999979$ &$0.14555794$ &$0.06757752$ &$0.03477145$ &$0.01139235$&$0.01789583$\\
14&$0.49977304$ &$0.14622041$ &$0.07041891$ &$0.03761547$ &$0.01139159$&$0.01788837$\\
15&$0.49977361$ &$0.14636154$ &$0.07274403$ &$0.04013040$ &$0.01139057$&$0.01788111$\\
$\cdots$&$\cdots$&$\cdots$&$\cdots$&$\cdots$&$\cdots$&$\cdots$\\
$\infty$&$0.50000000$ &$0.12500000$ &$0.05479177$ &$0.03097495$ &$0.01138938$&$0.01787406$\\
\hline
\end{tabular}
\end{center}
Note that there exists $\lim_{b\rightarrow\infty}\mu(a,b)$, and $\mu(a,b)\rightarrow 0$ uniformly in $b$, as $a\rightarrow\infty$.
Thus, the table above and some standard evaluations give the following
\begin{lem} Let $a,b,c\in\mathbb{N}$. Then
\begin{eqnarray*}
\mu(a,b)+\mu(b,c)\leq \mu(1,1)+\mu(1,9)<0.76.\quad\square
\end{eqnarray*}
\label{lema.1}
\end{lem}
\begin{lem}
There exists an absolute constant $c>0$ such that for all $\pp\in\mathbb{C}$, $\Re\pp\geq 1$, and all $a\in\mathbb{N}$, on has
$\big{|}\frac{\pp^{a}-1}{\pp-1}\big{|}>c$.
\label{lema.2}
\end{lem}
\proof Consider a contour, consisting of the segment $[1-iT,1+iT]$, and a semicircle $1+Te^{i\phi}$, $-\frac{\pi}{2}\leq\phi\leq\frac{\pi}{2}$.
For sufficiently large $T$, $\frac{\pp^{a}-1}{\pp-1}$ will be large on the semicircle. Moreover, this function never vanishes inside or on the
contour. Thus, from the maximum-minimum principle, its minimal absolute value is obtained on the segment $[1-iT,1+iT]$. Thus, let
$\pp=\frac{1}{\cos \psi}e^{i\psi}$, $-\frac{\pi}{2}<\psi<\frac{\pi}{2}$. Without loss of generality, let $\psi\geq 0$.
Consider the case $\frac{\pi}{2a}\leq\psi<\frac{\pi}{2}$. Then
\begin{eqnarray*}
\Big{|}\frac{\pp^{a}-1}{\pp-1}\Big{|}^{2}=\frac{\frac{1}{\cos^{2a}\psi}-\frac{2\cos a\psi}{\cos^a\psi}+1}
{\frac{1}{\cos^2\psi}-1}\geq \frac{\frac{1}{\cos^{2a}\psi}-\frac{2}{\cos^a\psi}+1}
{\frac{1}{\cos^2\psi}-1}=\frac{(\rho^a-1)^2}{\rho^2-1}:=Y(\rho),\quad \rho=\frac{1}{\cos\psi}.
\end{eqnarray*}
The function $Y(\rho)$ is an increasing function in $\rho$ for $\rho\geq 1$. It is obvious that
we may consider a case of $a$ sufficiently large. Thus,
\begin{eqnarray*}
\Big{|}\frac{\pp^{a}-1}{\pp-1}\Big{|}^{2}&\geq& Y\Big{(}\frac{1}{\cos\frac{\pi}{2a}}\Big{)}=\frac{\big{(}\frac{1}{\cos^a\frac{\pi}{2a}}-1\big{)}^2}{\tan^2\frac{\pi}{2a}}\\
&=&\frac{\Big{(}(1+\frac{\pi^2}{8a^2}+\frac{\mathcal{O}(1)}{a^3})^a-1\Big{)}^2}{\frac{\pi^2}{4a^2}+\frac{\mathcal{O}(1)}{a^3}}=
\frac{\frac{\pi^4}{64a^2}+\frac{\mathcal{O}(1)}{a^3}}{\frac{\pi^2}{4a^2}+\frac{\mathcal{O}(1)}{a^3}}=\frac{\pi^2}{16}+\frac{\mathcal{O}(1)}{a}.
\end{eqnarray*}
Let now $0\leq\psi<\frac{\pi}{2a}$. First, consider a function $\frac{1}{y}\log\cos (y\psi):=V(y)$. It is a decreasing
function for $0<y<\frac{\pi}{2\psi}$. Indeed, this is equivalent to the inequality
\begin{eqnarray*}
-\tan x\cdot x-\log\cos x<0, \text{ for }0< x<\frac{\pi}{2}.
\end{eqnarray*}
The function on the left is itself a decreasing function, with maximum value attained at $x=0$. Thus, $V(1)\geq V(a)$, which means $\cos a\psi\leq\cos^a\psi$, and this gives
\begin{eqnarray*}
\Big{|}\frac{\pp^{a}-1}{\pp-1}\Big{|}^{2}\geq\frac{\frac{1}{\cos^{2a}\psi}-1}{\frac{1}{\cos^2\psi}-1}\geq 1.\quad\square
\end{eqnarray*}
Therefore, Lemma \ref{lema.2} implies that the function $\pp^{-1}W^{-1}_{a}(\pp)$ is uniformly bounded:
\begin{eqnarray*}
\sup_{a\in\mathbb{N},|\pp-2|\leq 1}|\pp^{-1}W^{-1}_{a}(\pp)|=c_{0}<+\infty.
\end{eqnarray*}
This shows the validity of the following Lemma (apart from a numerical bound, which is the outcome of computer calculations).
\begin{lem} One has
\begin{eqnarray*}
\sup\limits_{|\pp-2|\leq 1,a\in\mathbb{N},|z-1|\leq 1}|\pp^{-1}W^{-1}_{a}(\pp)z-1|<1.29.\quad\square
\end{eqnarray*}
\label{lema.3}
\end{lem}
 \subsection{Numerical values for the moments}
\label{page14}
Unfortunately, Corollary \ref{cor1.3} is not very useful in finding exact decimal digits of $m_{2}$. In fact, the vector $(m_{1},m_{2},m_{3}...)$
is the solution of an (infinite)
system of linear equations, which encodes the functional equation (\ref{funk}) (see \cite{ga1}, Proposition 6). Namely, if we denote
$c_{L}=\sum_{n=1}^{\infty}\frac{1}{2^{n}n^{L}}=\text{Li}_{L}(\frac{1}{2})$,
we have a linear system for $m_{s}$ which describes the coefficients
$m_{s}$ uniquely:
\begin{eqnarray*}
m_{s}=\sum\limits_{L=0}^{\infty}(-1)^{L}c_{L+s}\binom{L+s-1}{s-1}m_{L},
\quad s\geq 1.
\end{eqnarray*}
Note that this system is not homogeneous ($m_{0}=1$). We truncate this matrix at sufficiently high order to obtain float values.
The accuracy of this calculation can be checked on the test value $m_{1}=0.5$.  This approach yields (for the matrix of order $325$):
\begin{eqnarray*}
m_{2}&=&0.2909264764293087363806977627391202900804371021955943665492_{+},\\
m_{3}&=&0.1863897146439631045710466441086804351206556532933915498238_{+}\\
m_{4}&=&0.1269922584074431352028922278802116388411851457617257181016_{+}.\\
\end{eqnarray*}
with all $58$ digits exact (note that $3m_{2}-2m_{3}=0.5$). In fact, the truncation of this matrix at an order $325$ gives
rather accurate values for
$m_{L}$ for $1\leq L\leq 125$, well in correspondence with an asymptotic formula \cite{ga3}
\begin{eqnarray}
m_{L}=\sqrt[4]{4\pi^{2}\log2}\cdot c_{0}\cdot L^{1/4}{\sf C}^{\sqrt{L}}+O({\sf
C}^{\sqrt{L}}L^{-1/4}),
\label{asympt}
\end{eqnarray}
where $c_{0}=\int_{0}^{1}2^{x}(1-F(x))\d x=1.030199563382_{+}$. So obtained numerical values for higher moments
tend to deviate from this expression rather quickly.\\

Kinney \cite{kinney} proved that the Hausdorff dimension of growth points of $?(x)$ is equal to
\begin{eqnarray*}
\alpha=\frac{1}{2}\Big{(}\int\limits_{0}^{1}\log_{2}(1+x)\d ?(x)\Big{)}^{-1}.
\end{eqnarray*}
Lagarias \cite{lagarias2} gives the following estimates: $0.8746<\alpha<0.8749$. Tichy and Uitz \cite{tichy_uitz} calculated $\alpha\approx 0.875$.
Paradis et al. \cite{paradis1} give the value $\alpha\approx 0.874832$.
We have (note that $?(1-x)+?(x)=1$):
\begin{eqnarray*}
A:=\int\limits_{0}^{1}\log(1+x)\d ?(x)=\int\limits_{0}^{1}\log\Big{(}1-\frac{1-x}{2}\Big{)}\d ?(x)+\int\limits_{0}^{1}\log 2\d?(x)=\\
-\sum\limits_{L=1}^{\infty}\frac{1}{L\cdot2^{L}}\int\limits_{0}^{1}(1-x)^{L}\d ?(x)+\log 2=-\sum\limits_{L=1}^{\infty}\frac{m_{L}}{L\cdot2^{L}}+\log 2.
\end{eqnarray*}
Thus, we are able to present much more precise result:
\begin{eqnarray*}
\alpha=\frac{\log 2}{2A}=0.874716305108211142215152904219159757...
\end{eqnarray*}
with all $36$ digits exact. The author of this paper have contacted the authors of \cite{paradis1} inquiring about the error bound for the numerical value
of $\alpha$ they obtained. It appears that for this purpose 10 generations of (\ref{cw}) were used. The authors of \cite{paradis1} were very kind in agreeing
to perform the same calculations with more generations. Thus, if one uses $18$ generations, this gives $0.874716 <\alpha<0.874719$.\\
\indent Additionally, the constant $c_{0}$ in (\ref{asympt}) is given by
\begin{eqnarray*}
c_{0}=\int\limits_{0}^{1}2^{x}(1-F(x))\d x=\frac{\mathfrak{m}(\log 2)}{2\log 2}=\frac{1}{2}\sum\limits_{L=0}^{\infty}\frac{m_{L}}{L!}(\log 2)^{L-1}.
\end{eqnarray*}
This series is fast convergent, and we obtain
\begin{eqnarray*}
c_{0}=1.03019956338269462315600411256447867669415885918240...
\end{eqnarray*}
\subsection{Rational functions $\mathbf{H}_{n}(z)$ }

The following is MAPLE code to compute rational functions
$\mathbf{H}_{n}(z)$={\tt h[n]}
and coefficients $\mathbf{H}_{n}'(0)$={\tt alpha[n]} for $0\leq n\leq 50$.\\
\line(1,0){480}
\begin{verbatim}
> restart;
> with(LinearAlgebra):
> U:=50:
> h[0]:=1/(2-z):
> for n from 1 to U do
>    j[n]:=1/2*simplify(
>        add( unapply(diff(h[n-j],z$j),z)(2*z)*2/j!*(z^(j)),j=1..n)+
>        add( unapply(diff(h[n-j-1],z$j),z)(2*z)*1/j!*(z^(j)),j=1..n-1)+
>             unapply(h[n-1],z)(2*z)    ):
>    k[n]:=simplify((z-1)^(n+1)*(unapply(j[n],z)(z)-
>                                unapply(j[n],z)(1/z)/z^2)):
>    M[n,1]:=Matrix(n,n):M[n,2]:=Matrix(n,n): M[n,3]:=Matrix(n,n):
>          for tx from 1 to n do for ty from tx to n do
>                    M[n,1][ty,tx]:=binomial(n-tx,n-ty)
>                         end do: end do:
>          for tx from 1 to n do M[n,2][tx,tx]:=2^(n-tx) end do:
>          for tx from 1 to n do M[n,3][tx,n+1-tx]:=2^(tx-1) end do:
>    Y[n]:=M[n,1]-1/2^(n+1)*M[n,2]+(-1)^(n+1)/2^(n+1)*M[n,3]:
>    A[n]:=Matrix(n,1):
>          for tt from 1 to n do A[n][tt,1]:=coeff(k[n],z,n-tt) end do:
>    B[n]:=MatrixMatrixMultiply(MatrixInverse(Y[n]),A[n]):
>    h[n]:=add(z^(n-s)*B[n][s,1](s,1),s=1..n)/(z-2)^(n+1):
>            end do:
>
> for n from 0 to U do alpha[n]:=unapply(diff(h[n],z$1),z)(0) end do;
\end{verbatim}
\line(1,0){480}\\
\indent It causes no complications to compute {\tt h[n]} on a standard home
computer for $0\leq n\leq 60$, though the computations heavily increase in
difficulty for $n>60$.

\subsection{Rational functions $\mathbf{Q}_{n}(z)$}

This program computes $\mathbf{Q}_{n}(z)=${\tt q[n]} and the values\\
$\mathbf{Q}'_{n}(-1)=${\tt beta[n]} for $0\leq n\leq 50$.\\
\line(1,0){480}
\begin{verbatim}
> restart;
>q[0]:=-1/(2*z);
>N:=50:
>q[1]:=simplify(1/2*unapply(q[0],z)(-1)*(1-1/z^2)):
> for n from 1 to N do
>   q[n]:=1/2*simplify(
>       add(unapply(diff(q[n-j-1],z$j),z)(-1)/j!*(z^(j)-1/z^(j+2)),j=1..n-1)+
>       unapply(q[n-1],z)(-1)*(1-1/z^2)
>                      ):
                end do:
> for w from 0 to N do beta[w]:=unapply(diff(q[w],z$1),z)(-1) end do;
\end{verbatim}
\line(1,0){480}\\

\par\bigskip

\noindent Vilnius University, The Department of Mathematics and Informatics, Naugarduko 24, Vilnius, Lithuania\\

\noindent Max-Planck-Institut f\"{u}r Mathematik, Vivatsgasse 7, 53111 Bonn, Germany\\

\noindent\tt{giedrius.alkauskas@gmail.com}
\smallskip
\end{document}